\documentclass[11pt]{article}
\usepackage{graphicx}
\usepackage[english]{babel}
\usepackage{theorem}

\advance\textwidth by 4cm
\advance\oddsidemargin by -2cm
\advance\evensidemargin by -2cm

\theorembodyfont{\rmfamily}

\advance\textheight by 2cm
\advance\topmargin by -1cm

\newcommand{\ol}{\overline}

\newcommand{\cvd}{$\| \bigskip$}

\newcommand{\T}[1]{\texttt{#1}}

\newenvironment{GAPproc}[1]{\noindent\texttt{#1}\par\smallskip\noindent}{\bigskip}

\newtheorem{theorem}{Theorem}[section]
\newtheorem{prop}{Proposition}
\newcommand{\bp}{\begin{prop}}
\newcommand{\ep}{\end{prop}}
\newtheorem{defi}[prop]{Definition}

\newtheorem{question}{Question}
\newcommand{\bd}{\begin{defi}}
\newcommand{\ed}{\end{defi}}
\newcommand{\bde}{\begin{description}}
\newcommand{\ede}{\end{description}}
\newcommand{\bi}{\begin{itemize}}
\newcommand{\ei}{\end{itemize}}
\newtheorem{conj}{Conjecture}
\newcommand{\bc}{\begin{conj}}
\newcommand{\ec}{\end{conj}}
\newcommand{\be}{\begin{enumerate}}
\newcommand{\ee}{\end{enumerate}}
\newcommand{\bt}{\begin{theorem}}
\newcommand{\et}{\end{theorem}}
\newcommand{\bdes}{\begin{description}}
\newcommand{\edes}{\end{description}}
\newcommand{\ba}{\begin{array}}
\newcommand{\ea}{\end{array}}
\newcommand{\bit}{\begin{itemize}}
\newcommand{\eit}{\end{itemize}}

\title{Towards an implementation of the B-H algorithm
\\for recognizing the unknot}

\author{J. S. Birman, P. Boldi, M. Rampichini, S. Vigna}
\begin{document}

\maketitle

\begin{abstract}
In the manuscript \cite{bh} the first author and Michael Hirsch presented a then-new 
algorithm for recognizing the unknot. The first part of
the algorithm required the systematic enumeration of all discs which support a `braid
foliation' and are embeddable in 3-space. The boundaries of these `foliated embeddable
discs' (FED's) are the collection of all closed braid representatives of the unknot, up
to conjugacy, and the second part of the algorithm produces a word in the
generators of the braid group which represents the boundary of the previously listed 
FED's. The third part tests whether a given closed braid is conjugate to the boundary of
a FED on the list. 

In this paper we describe implementations of the first and second parts of the 
algorithm. We also give some of the data
which we obtained. The data suggests
that FED's have unexplored and interesting structure. Open questions are interspersed
throughout the manuscript.  

The third part of the algorithm was studied in
\cite{bkl:1} and \cite{bkl:2}, and implemented by S.J. Lee \cite{l}. At this writing his
algorithm is polynomial for $n\leq 4$ and exponential for $n\geq 5$. 
\end{abstract}

\section{Introduction}
\label{section:introduction}

\subsection{Background}

The subject of this paper is the question: given a knot $K$, can we decide
whether $K$ is the unknot?  The problem was solved affirmatively by W.
Haken in a groundbreaking paper published in 1961 \cite{ha}.  However, showing 
that an algorithm exists does not mean that there is necessarily an algorithm
which will be useful in practice, even for the simplest examples.  Thus in 1993
(over 30 years after Haken did his work) when Hoste, Thistlethwaite and Weeks 
tabulated the 1,701,936 prime knots  with $\leq$16 crossings \cite{htw} they 
had all the tools of the trade available to them, but 
used a ragbag of diagrammatic techniques to eliminate unwanted
appearances of the unknot. In that regard it should be noted that knot
diagrams with at most  16 crossings do not even begin to exhibit the pathology 
which one knows exists in the general case. For example, see \cite{g} for some 
examples which show why the  diagrammatic approach was abandoned in the 1930's. 
(On the other hand, see \cite{hl} for a recent proof that an upper bound exists 
for the number of Reidemeister moves which must be tested to be sure that a 
knot diagram with a given number of crossings is not the unknot.)

Haken's work begins by constructing a triangulation of the complement of $K$. 
He then applies the theory of normal surfaces, due to Kneser \cite{kn}, who 
showed that any surface $F$ of minimal genus with boundary $K$ can be assumed 
to be in a special position in which it intersects each tetrahedron $T_i$ in 
the triangulation in an especially nice way, namely as a set of parallel sheets,
each sheet being a polygonal disc whose boundary has 3 or 4 edges in the 
faces of $\partial T_i$.  The polygonal discs are used to set up a system of 
linear equations. Solving the system allows one to decide whether, in fact, the
solution set includes a normal surface which is a disc.

A very different approach to the unknot recognition problem was discovered by
the first author and M. Hirsch, who developed in \cite{bh} the algorithm which
is the subject of investigation in this paper.  The basic idea behind the B-H
algorithm rests in the braid foliation techniques of Birman and Menasco (see 
\cite{bm5} and \cite{bf}). Braid foliations allow one to generate, in a 
systematic manner, a list of all of the foliated embedded discs whose 
boundaries are closed braid representatives of the unknot. The list is ordered 
by a complexity function which depends on properties of the foliated discs. 
One then compares a given example $K$ with the examples on the list in order to 
decide whether $K$ is the unknot.

In this paper we will give a computer implementation of certain parts of the
algorithm in \cite{bh}, namely the problems of enumerating the foliated 
embeddable discs and finding the braid words which describe their boundaries. 
We note that for  braid index $2$ the problem is trivial. For braid index $3$ 
the unknot recognition problem was solved in \cite{mp}, where it was proved 
that there are precisely 3 conjugacy classes of closed 3-braid representatives 
of the unknot. For $n=4$ the question is much harder, because of the example 
in \cite{m:4b} and the others which are presented here. See \cite{f} for a 
proof that there are infinitely many distinct conjugacy classes of 4-braid
representatives of the unknot. We were able to obtain non-trivial data for
braid index $n=4$. The examples which
we found have braid word descriptions with $\leq$11 band generators, however 
there may well be shorter braid words for the same examples. We also give a 
small amount of scattered data for higher braid index. We note that for $n=4$ 
polynomial-time algorithms exist for the solution to the conjugacy problem and 
the shortest word problem \cite{kkl}, which could easily be integrated with our 
work, however we did not make a systematic attempt to do that.

The data which we obtained is given in Section~\ref{sec:data} of this paper.

We conjecture that a practical polynomial time algorithm exists which will
solve the unknot recognition problem in the special case of knots of braid 
index 4. Theorem 4.3 of the review article \cite{bf} (which gives a new proof 
of the main result in \cite{bm5}) would surely play an important role in any 
such solution, as would the polynomial-time algorithms of \cite{kkl}. The chief obstacle, as we see it, is to find an efficient way to enumerate all the 
foliated embeddable discs with $N$ negative vertices and $4+N$ positive 
vertices.  We believe that when the structure of these foliated embeddable 
discs is better understood, this problem will be solved.

\subsection{A review of braid foliations}

In this section we briefly review the main results of \cite{bh}. 
A good reference for a survey on braid foliations is \cite{bf}.
After completing our review of the results which we need from \cite{bh} and \cite{bf} we
explain in a precise way what we do in this paper.

The underlying plan is the following: the unknot is the unique knot which bounds a disc
embedded  in $\mathbf{R}^3$. All discs embedded in $\mathbf{R}^3$ can be isotoped in
such a way that the boundary is a closed braid relative to the $z$-axis. All these
embeddings can be described by a finite set of combinatorial data, and they can be
listed in order of increasing complexity. To each disc we will show how to associate the
braid whose closure is the  boundary of the disc. 
So if we want to know whether a given knot $K$ is the unknot, we first 
represent it as a closed braid $\hat\beta$, using our preferred algorithm (see 
\cite{m:tk,v,y}). Then compare this braid with our list of braids which are
the boundary 
of an embedded disc, looking for a braid $\gamma$ which is in the same conjugacy class
as the given braid  $\beta$. To check conjugacy, use for instance the algorithm in
\cite{bkl:1} and \cite{bkl:2}, as implemented by S.J. Lee \cite{l}.  In
\cite{bh} an upper bound is given for the complexity of the disc to look for in the
list, so the process is finite.  The problem of improving that upper bound will not be
considered here.

To implement our algorithms we have used [GAP 99] The GAP Group, GAP --- Groups, Algorithms, and Programming,
Version 4.1; Aachen, St Andrews, 1999.

\ 

Let us fix the $z$-axis $A$ of $\mathbf{R}^3$ as braid axis, and the standard fibration by half-planes 
of the complement of $A$ by half-planes  
$$H_\theta=\{(\rho\cos\theta,\rho\sin\theta,z) \, | \, 
\rho>0, z\in\mathbf{R}\}$$ 
\bt[cf Theorem 2.1 of \cite{bh}] A disc $D$ embedded in the standard fibration
of 
$\mathbf{R}^3\setminus A$, with boundary a closed braid
$\hat\beta$, $\beta\in B_n$, can always be put in general position so that the induced
singular foliation on $D$ has the following properties:
\be
\item All intersections of the disc $D$ with the axis $A$ are transversal.
These intersections consist of $P$ positive and $N$ negative points, where the sign
is positive if the orientation of $A$ agrees with that of $D$, otherwise negative. 
We call these
$P+N$ intersection points `vertices'. The braid index is $n=P-N$.
\item The disc $D$ intersects almost all half-planes $H_\theta$ transversally, in
what we call regular leaves.
\item The foliation in a neighborhood of $\partial D$ is transverse.
\item The foliation in a small circular neighborhood of each vertex is radial.
\item There are a finite number $P+N-1$ of singular half-planes $H_{\theta_i}$ to
which $D$ is tangent  in one point, which is a non-degenerate saddle. A saddle
together with its four leaves (branches) is called a singular leaf.
\item The branches of a saddle can be of the following two types: type $a$: a
simple arc with one endpoint on $\hat\beta$ and the other on $A$; type $b$: a
simple arc with both endpoints on $A$.
\item The saddles are restricted to the following types:
\begin{enumerate}
\item [] $aa$-saddles: singularities between two $a$-arcs;
\item [] $ab$-saddles: singularities between an $a$-arc and a $b$-arc;
\item [] $bb$-saddles: singularities between two $b$-arcs.
\end{enumerate}
Each saddle can be either positive or negative, according as  the orientation of the
disc and the tangent half-plane at the saddle point agree or disagree. 
\item The vertices are cyclically ordered along $A$, and the saddles are
cyclically ordered around $A$.
\ee
\et

\noindent \textbf{The code for a foliated disc:} It is not explicitly
explained  in \cite{bh} how to represent an embeddable disc on a computer:
we will encode the description of an embedded disc $D$ in terms of its vertices 
and saddles as follows:
\be
\item [] {\bf The vertex string $V$} is the list of $(P,N)$ vertices; each
positive vertex will be denoted by an integer number 
$k\in\{1,2,\ldots P\}$; each negative vertex will be denoted by 
a pair of integer numbers $k.j$, where $k$ is the number of the immediately 
preceding positive vertex in $V$ ($k$ might be 0) and $j$ is the 
ordinal number of the negative vertex in the subset of 
negative vertices between $k$ and $k+1$;
\item [] {\bf The ordered list of saddles,} in which each saddle will be denoted
by:   the list of vertices involved in the saddle: 
two positive for an $aa$-saddle, two positive and one negative for an $ab$-saddle,
two positive and two negative for a $bb$-saddle; and the sign, $\pm1$. 
\ee

For a fixed braid index $n$ we can assign to a disc $D$ its \textbf{complexity} $(P,N)$,
that is the 
number of positive and negative vertices, with $P-N=n$, and list all embeddable discs in order 
of increasing complexity. This is possible thanks to the following theorems:

\bt[cf Theorem 2.2 of \cite{bh}]
The combinatorial data for an embeddable disc $D$, i.e. the cyclically ordered
list of vertices, with  their signs, and the cyclically ordered list of saddles
with their signs, determine the embedding in
$\mathbf{R}^3$, uniquely up to foliation-preserving isotopy. They also determine the embedding of 
the boundary of the disc as a closed braid.
\et

To see an example of the singular foliation on an embeddable disc, look ahead to
Figure \ref{thedisc}:  We show there a disc $D$ in which we have drawn
all vertices and saddles, with all singular leaves.  
The complement of the
singular leaves in the disc is the disjoint union of open discs,  some of them
bounded by: an arc of
$\partial D$, a positive vertex  and some singular leaves (there are 18 like
this in Figure \ref{thedisc}); the others (maybe zero) bounded by one positive 
and one negative vertex, and some singular leaves (there are 13 like this in Figure \ref{thedisc}). Discs of the first type 
can be foliated by regular leaves which are $a$-arcs, with endpoints on the positive vertex 
and $\partial D$. Discs of the second type can be foliated by regular leaves which are $b$-arcs, 
with endpoints on the positive and the negative vertex.

\bt[cf Theorem 3.4 of \cite{bh}] From the set of combinatorial data describing an embeddable disc
$D$ we extract a unique extended boundary word, which represents a braid whose closure is the link
consisting of the boundary of $D$ and the $N$ unlinked small circles bounding small disc neighborhoods
of the $N$ negative vertices of $D$.
\et

In section \ref{boundaryword} we will give an explicit implemented algorithm to get the boundary braid 
from the extended boundary braid.

For a given set of combinatorial data as described before we have to test embeddability. 
For this purpose, let us give some necessary definitions:
\bd[cf \cite{bh}] A {\bf b-arc} is a regular leaf $b(i,j.h)$ which is a
simple arc  connecting a positive vertex $i$ to a negative vertex $j.h$. 
A {\bf generalized  or  gb-arc} is either a $b$-arc or the part of a
singular leaf of an $aa$-saddle connecting the two positive vertices and passing
through  the saddle point.
\ed

Since the foliation around each vertex is the standard radial foliation, we can
distinguish leaves around each vertex by means of their angle $\theta$. If there
exists at some $\theta$ a $b$-arc $b(i,j.h)$, then there exists a maximal
interval
$(\theta_l,\theta_k)$ in which all regular leaves around $i$  and $j.h$ are
$b(i,j.h)$. In this case we say that the $b$-arc exists in
$(\theta_l,\theta_k)$. Say that the $gb$-arc $gb(i,j)$ exists in
$(\theta_{k-1},\theta_k)$ if there is an $aa$-saddle $(i,j)$ occuring at
$\theta_k$.
For instance, looking at Figure \ref{newdisc} see $b(3,0.1)$ in $(\theta_{14},\theta_7)$ 
(leaves go clockwise around a negative vertex) and $gb(4,8)$ in $(\theta_5,\theta_6)$.

\bt[Theorem 3.5 of \cite{bh}]
\label{embeddability theorem}
A disc $D$, given in terms of combinatorial data, is embeddable if and only if it satisfies the 
following three conditions:
\be
\item The saddles about each positive (respectively negative) vertex  are in
counterclockwise (respectively clockwise) order;
\item  The vertices attached to a positive (respectively negative) saddle are in
counterclockwise (respectively clockwise) order;
\item The endpoints of a $(g)b$-arc in $(\theta_{k-1},\theta_k)$ never separate
the endpoints of a $b$-arc in the same interval.
\ee
\et

In what follows the expression `a disc with $(P,N)$ vertices' will mean
a disc (or its combinatorial description) with $P$ positive and $N$ negative vertices.

The fundamental tool for listing all embeddable discs in \cite{bh} is the \textbf{insertion of an $ab$-tile}:
given an embeddable disc $D$ with $(P,N)$ vertices, Birman and Hirsch describe how to get from $D$ a new disc
$D'$ with $(P,N+1)$ vertices and one more saddle of type $ab$.

\bt[cf Theorem 4.1 of \cite{bh}] Each embeddable disc $D$ with $(P,N)$ vertices can be constructed by
starting from an embeddable disc $D_0$ with $(P,0)$ vertices, adding $N$ $ab$-tiles one at a time. At each 
stage the new negative vertex and the new $ab$-saddle
are inserted into the order of the older vertices and saddles, in such a way
that the new disc is embeddable.
\et

\bt[Theorem 4.2 of \cite{bh}]\label{enumeration bh}
All possible embeddable discs of fixed braid index $n$ may be enumerated
in order of increasing $(P,N)$, with $P-N=n$, by the following (not necessarily efficient) procedure:
\bi
\item enumerate all possible discs with $(n,0)$ vertices, testing each for embeddability;
discard non embeddable ones;
\item enumerate all discs with $(n+j,0)$ vertices, testing each for embeddability;
discard non embeddable ones; then add $j$ $ab$-tiles in all possible ways, testing each obtained disc 
for embeddability; discard non embeddable ones; get all embeddable discs with $(n+j,j)$ vertices.
\ei
\et

We will call discs with zero negative vertices \textbf{positive discs}.
Since we are interested in conjugacy classes of braids, we have to notice that this list has duplicates, 
because non isotopic embeddable discs may have the same boundary braid, or different
boundary braids in the same conjugacy class.

In what follows we will show:
\bit
\item how our set of combinatorial data describes an (embeddable) disc;
\item how to translate many properties of the combinatorial foliation in 
braid words in the band generators;

\item how to enumerate all positive embeddable discs;
\item how to reduce a lot of the redundancy in the resulting list;
\item how to insert $ab$-tiles and test embeddability with a computer program;
\item how to find the boundary word of an embeddable disc;
\item how to fill the list via a different way: no more by insertion of $ab$-tiles and embeddability test,
but via enumeration of sequences of half-planes completely describing 
embeddable discs.
\eit
The main reason for wanting to see the data produced by an algorithm is if it
suggests new structure which will then lead to a better algorithm. For this reason we
will intersperse ``questions'' throughout the manuscript, as they occur to us in the
context of our work. Some of them have immediate answers, but most relate to open
problems.

\section{The disc, its code and\\the associated braid word}
\label{drawthedisc}

We will use the band presentation of the braid group which is given in \cite{bkl:1}. The
generators are the $n\choose2$ braids $\{(i,j), n\geq i>j\geq1\}$, such that the $i^{th}$
strand crosses over the
$j^{th}$ strand, with all the other strands left unchanged under these two.
The relations are of two types:
$$(i,j)(k,l)=(k,l)(i,j), \mbox{ when } (i-k)(i-l)(j-k)(j-l)>0,$$
(the condition means that the two pairs of indices are non interlocking), and
$$(i,j)(j,k)=(i,k)(i,j)=(j,k)(i,k), \mbox{ when } n\geq i>j>k\geq1.$$
We will denote the inverse of a generator $(i,j)$ by $\ol{(i,j)}$.

\ 

Suppose we are given a list of combinatorial data for an embeddable disc. In this section we explain
how to draw the singular foliation of the disc, and how to associate to $D$ its extended boundary word. 
Let us consider the following example: 

$$D=\{[[4.1,6,7],1],[[9,10],1],[[8.2,10,11],1],[[4,8.1,11],1],[[4,8],1],$$
$$[[0.1,3,7],-1],[[0.1,7,8],-1],[[0.2,2,4.1,6],-1],[[0.2,6,7],-1],[[5,6],1],$$
$$[[1,5],1],[[0.2,2,4.1,7],1],[[0.1,3,8],1],[[4,8.1,11],-1],[[8.2,10,11],-1]\}.$$

\bd[cf \cite{bh}]
For an embeddable disc $D$ with $(P,N)$ vertices, its {\bf extended
boundary word}  is the braid word in the band generators
of $B_P$ with each letter given by the pair of positive vertices and the sign of 
the corresponding saddle of $D$. 
\ed
In our example, 
 
$$EW(D)=(7,6)(10,9)(11,10)(11,4)(8,4)\ol{(7,3)}\ol{(8,7)}\ol{(6,2)}\cdot$$
$$\cdot\ol{(7,6)}(6,5)(5,1)(7,2)(8,3)\ol{(11,4)}\ol{(11,10)}.$$

\bp\label{extended word}
If a word $W$ of a braid in $B_P$ is the extended boundary word of an embeddable disc
with $(P,N)$ vertices, then it has the following properties:
\be
\item The length of $W$ (in band generators) is $P+N-1$;
\item The induced permutation $\rho(W)\in S_P$ is a product of one $(P-N)$-cycle
and $N$ 1-cycles.
\ee
\ep

\noindent\textbf{Proof:} To each saddle of the disc corresponds one letter of $W$, 
in the same order, 
to each positive vertex corresponds one index, and to each negative 
vertex corresponds one
single unknotted strand of the closed braid, unlinked from the rest.
The first condition comes from the Euler characteristic of the 
foliated disc. \cvd

\begin{figure}
\begin{center}\includegraphics{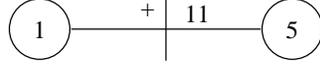}\end{center}
\caption{The first saddle.}
\label{firstsaddle}
\end{figure}

\begin{figure}
\begin{center}\includegraphics{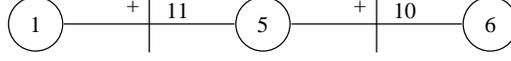}\end{center}
\caption{The first two saddles.}
\label{twosaddles}
\end{figure}

In our example $EW\in B_{11}$, so that $P=11$; its length is 15, therefore 
$N=5$ and $P-N=6$; 
the associated permutation is 
$$\rho(EW)=(1,5,6,8,10,9)(2)(3)(4)(7)(11).$$
Now we can explain how to draw  
leaves of the foliation 
on the disc. 
Consider the first index of the $(P-N)$-cycle, say $j_1$: look for the first letter 
of $EW(D)$ in which that index occurs: say $(j_1,j_2)^{\pm1}$ (or $(j_2,j_1)^{\pm1}$). 
Draw the two positive vertices and the singular leaves joining them, and label the vertices
with their numbers and the saddle with its sign and number (see Figure \ref{firstsaddle}). 
Then look for the (cyclically) next letter containing $j_2$: 
if it is a different one, then draw it attached to the previous one 
(see Figure \ref{twosaddles}), if it is the same letter, then procede to the next
letter containing $j_1$.  We have to run twice through each $aa$-saddle.

\begin{figure} 
\begin{center}\includegraphics{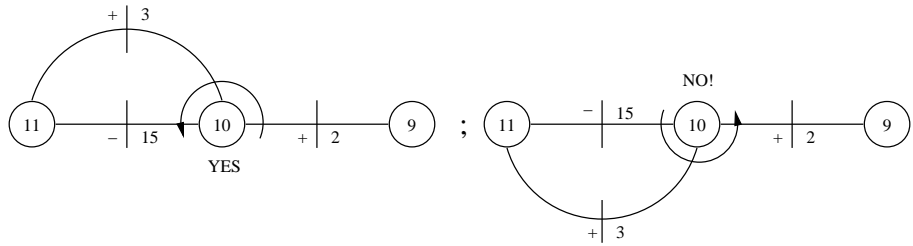}\end{center}
\caption{Respect order about each vertex.}
\label{respectorder}
\end{figure}

To respect embeddability, when a positive vertex has two or more saddles attached, 
and we have to draw another one attached to it, we must put it in the right 
(counterclockwise) cyclic order about the vertex (see Figure \ref{respectorder}).

\begin{figure}
\begin{center}\includegraphics{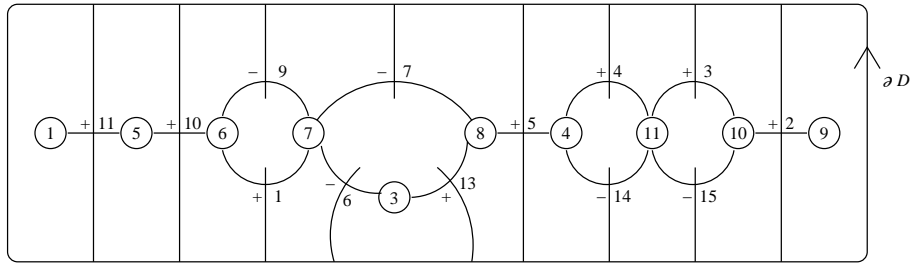}\end{center}
\caption{The disc with all its `external' saddles.}
\label{discwithholes}
\end{figure}

\begin{figure}
\begin{center}\includegraphics{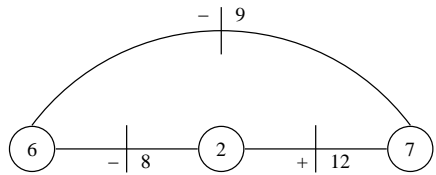}\end{center}
\caption{The cycle of saddles about another positive vertex.}
\label{cycleonevertex}
\end{figure}

Procede as indicated above until you get back to $j_1$ and the $(P-N)$-cycle is
completed. These singular leaves divide the disk into an outer part,
which is connected to the boundary, and some inner parts, in which
the negative vertices lie along with possibly those positive vertices which at this
stage  have not yet been drawn (see Figure \ref{discwithholes}). 
They are vertices occurring in the 1-cycles of the
permutation. For each of them, draw its cycle of saddles in a similar
way (see Figure \ref{cycleonevertex}), then attach it inside the appropriate inner region. 

\begin{figure}
\begin{center}\includegraphics{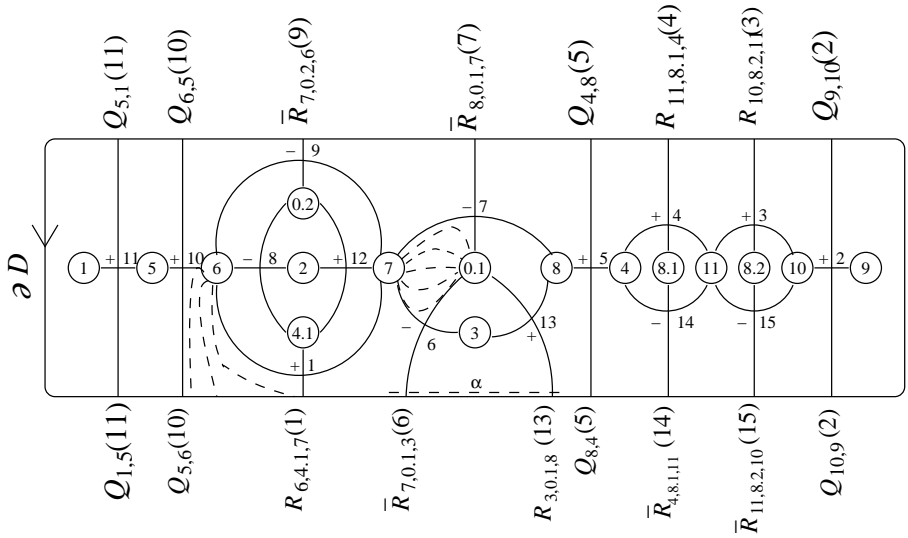}\end{center}
\caption{The foliated disc and the code for its boundary.}
\label{thedisc}
\end{figure}

It remains to add the negative vertices. Look at the code of the
disc, and put each negative vertex in its place, joining it to the
(already existing) saddles by the other singular leaves (see Figure \ref{thedisc}).

\

\noindent \textbf{Remark:} Notice that this process of drawing a disc can be performed
almost entirely from the extended boundary word. Only at the end do we need to know
the exact position of the negative vertices in the vertex string. But their number $N$
and their topological position inside the inner regions of the foliation are already
specified by $EW$. In some cases it is also possible to decide their position in the
vertex string simply by reading their relative order about the saddles to which they
are attached. 

\begin{figure}
\begin{center}\includegraphics{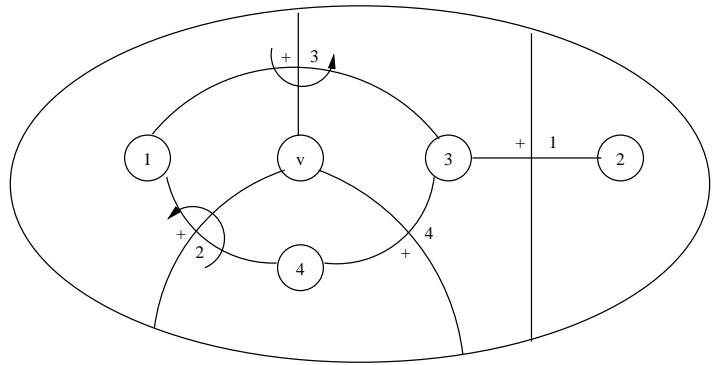}\end{center}
\caption{A non embeddable disc.}
\label{wrongdisc}
\end{figure}

\begin{question}
It's natural to ask whether every word which satisfies the conditions of Proposition 
\ref{extended word} actually can be realized by an embeddable foliated disc?
The answer is `no'. For example $W=(3,2)(4,1)(3,1)(4,3)$ satisfies the
conditions of Proposition \ref{extended word} but the corresponding disc is not
embeddable.  To see this consult Figure \ref{wrongdisc}. For saddle 2 we must have $4<v<1$,
but for saddle 3 we require $1<v<3$, which is impossible. 
\end{question}

If $EW=(6,5)(5,4)\ol{(4,2)}(3,1)(5,3)(4,2)$, we can draw all the disc and see that the negative
vertex must be either $0.1$ or $5.1$. But neither of them is embeddable: they both pass
the first part of embeddability test, but not the $b$-arc test: in the first case, in
$(4,5)$ we have $gb(3,5)$ interlocking with $b(4,0.1)$; in the second case we have
for instance in $(3,4)$ $gb(3,5)$ interlocking with $b(4,5.1)$.

\begin{question}
Are there other necessary conditions on the extended boundary word EW? For example, we
know that each pair of positive vertices may be involved in at most two $ab$-saddles
one opposite to the other, or just in one saddle. Therefore each letter must appear at
most once, also its inverse may appear, but if so exactly once.  A more efficient
algorithm would clearly be possible if we knew a better set of conditions that would
allow us to rule out certain discs on the basis of the associated boundary words.
\end{question}

\section{Positive discs, good words\\and how to reduce redundancy}
In a disc without negative vertices, condition (3) of Theorem \ref{embeddability
theorem} is vacuous, because there are no negative vertices, and so there are no
$b$-arcs.  
This means that the embeddability test is considerably simpler for positive words
than for arbitrary words, because the difficult  part is the test for $b$-arcs.

\bp
The word $W$ associated to a positive disc with $(n,0)$ vertices has the 
following properties:
\be
\item Its length in the band generators is $n-1$;
\item The induced permutation $\rho(W)$ is an $n$-cycle.
\ee
Moreover any word with these properties corresponds to  the boundary word of an
embeddable positive disc. Also, this word encodes all information about the disc.
\ep


\noindent\textbf{Proof:} The first part is a corollary of Proposition \ref{extended word}.  
The second part is known, because the braids in question are in fact the `Stallings
braids', which form   a proper subset of the braids whose closure is the unknot 
(see eg \cite{m:eb}). For the last sentence, as we have seen in the preceding section,
the extended word  alone carries all information except the order of negative vertices
in the vertex string.  But a positive disc has no negative vertices. \cvd

\bd 
We call any word satisfying the two conditions of the preceding Proposition
a \textbf{good word}.
\ed


As we anticipated in the remark at the end of the previous section, a good word
completely determines the embedding of its disc. Indeed, a description of this
embedding is computed by our GAP procedure GenerateDiscBoundary($n,W$).

\medskip

To follow the program of Theorem \ref{enumeration bh} we first need to list all
the positive discs, that correspond to all the positive good words. But since we are
interested in conjugacy classes of braids, we can reduce a lot  of redundancy at this
stage by some easy conjugations. 
The following definition concerns all embeddable
discs (not only the positive ones).
\bd 
We say that two embeddable foliated discs are {\rm {\bf equivalent}} if they only
differ by a cyclic permutation of the names of the vertices or saddles.
\ed

It is clear that two equivalent discs are isotopic. That is, not only are their
boundary words equivalent as cyclic words, but in fact the entire singular foliation
is the same, up to a cyclic permutation of the  `names' of the vertices and saddles. 

\bp
The (extended) boundary words $W,W'$ of two equivalent discs $D,D'$ only differ 
by some of the following easy conjugations:
\be
\item conjugations by powers of $\delta$, where $\delta=(n,n-1)(n-1,n-2)\ldots(3,2)(2,1)$, and 
$\delta^{-1}(i,j)\delta=(i+1,j+1)$ (mod $(n,n)$) for all band generators $(i,j)$ (see
{\rm\cite{bkl:1}});
\item conjugations by initial or final subwords;
\ee
\ep
\noindent\textbf{Proof:} Conjugations of the first kind correspond to cycling the names of
vertices along the braid axis; conjugations of the second kind correspond to cycling
the names of saddles about the braid axis.
\cvd

\ 

We call them `easy conjugations' because it is very easy and inexpensive to
perform them on a computer (see the Appendix).

So for a given $n$ we will list all good words up to easy conjugations and
inversions: see Proposition \ref{inversion}.  We have written a function in GAP
called EnumeratePositiveGoodWords($P$) that enumerates one representative for
each orbit of the action of the group $G=S_P\times S_{P-1}$ on the set of
positive words of $B_P$ with length $P-1$, where $S_P$ acts on the $P$ indices
and $S_{P-1}$ acts on the letters of the word, discarding those which have not
the required permutation property.

Once we have all positive good words up to easy conjugations, we can list all
positive good words up to inversion by choosing
$1,2,\ldots\lfloor\frac{P-1}2\rfloor$ letters of each word to become negative.
This is done by another easy procedure, called EnumerateGoodWords$(P)$.
Remember that words with different exponent sum are surely non conjugate.

For instance, the result of EnumerateGoodWords(4) is a list of 32 good words of exponent sum 3 
or 1.
We will get 32 other good words of exponent sum $-3$ and $-1$ by inversion. Notice
that among these 64 good words, some represent the same braid, because it is possible
to apply some relations: 
for instance $(2, 1)\ol{(3, 1)}(4, 3)=(2, 1)(4, 1)\ol{(3, 1)}$, 
and both words appear in the list; but notice that
$(2, 1)(3, 1)(4, 3)\neq(2, 1)(4, 1)(3, 1)$, also both appearing in the list.
So we cannot reduce the list also by relations, both because relations are different on words with different
signs, and they are not easy to be performed by the computer.

S. J. Lee has reduced the 32 good words up to conjugation: there are three conjugacy classes
with exponent sum 3, with representatives  $(2,1)(3,1)(4,1)$, $(2,1)(3,1)(4,2)$ and
$(2,1)(4,2)(3,1)$.  There are four conjugacy classes with exponent sum 1, with
representatives $\ol{(2,1)}(3,1)(4,1)$,  $\ol{(2,1)}(4,1)(3,1)$,   $\ol{(2,1)}(3,1)(4,2)$
and $(2,1)(4,2)\ol{(3,1)}$.

\begin{question}  Are easy conjugations and defining
relations in the braid group  the only  possible moves between conjugate positive
EW's of the same length?  A better understanding of this issue would lead to a more
efficient method of listing the positive words which we need to test.  As will be seen,
any redundancies which can be eliminated at this stage of the algorithm will lead to
major savings at subsequent stages, enabling us to collect better data.
\end{question}


\section{How to insert new $ab$-tiles in a given disc}
\begin{figure}
\begin{center}\includegraphics{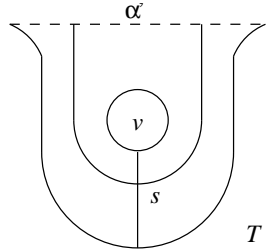}\end{center}
\caption{An $ab$-tile.}
\label{abtile}
\end{figure}

\begin{figure}
\begin{center}\includegraphics{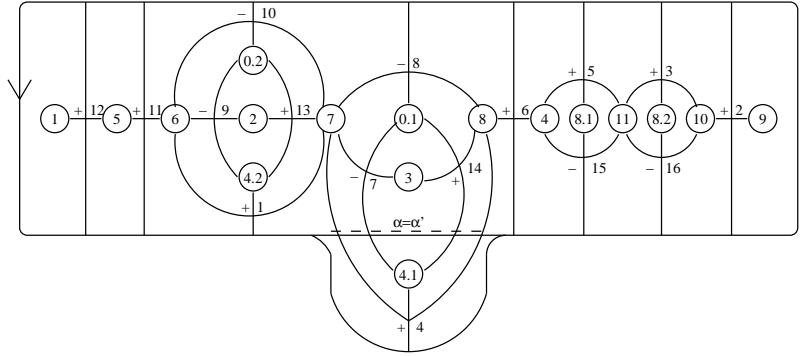}\end{center}
\caption{The foliated disc obtained by insertion of an $ab$-tile
 from the disc of Figure \ref{thedisc}.}
\label{newdisc}
\end{figure}

The process of inserting an $ab$-tile in a given foliated disc $D$ is explained well
in \cite{bh}. The idea is to take another small disc $T$ (see Figure \ref{abtile}),  with
a negative vertex $v$ and a saddle $s$ inside it, and a distinguished arc 
$\alpha'$ on its boundary $\partial T$, such that the four branches of $s$
end one in $v$, two in $\alpha'$ and the fourth in $\partial T\setminus\alpha'$.
Then choose an \textbf{insertion arc} $\alpha$ along the boundary $\partial D$ 
and attach $T$ along $\alpha$ by identifying $\alpha$ with $\alpha'$ (see Figure \ref{newdisc}), 
and continue branches of saddles of $D$ ending in $\alpha$ inside $T$ till they arrive
at $v$, and continue the two branches of $s$ ending in $\alpha'$ inside $D$
till they arrive at two specified positive vertices of $D$.
Get a new foliated disc $D'$ in such a way that this is still embeddable.

We know the conditions for embeddability from Theorem \ref{embeddability theorem},
so we can impose conditions on the choice of the insertion arc $\alpha$, on
the position of $v$ in the new vertex string of $D'$, and the position of
$s$ in the new saddle list of $D'$, so that at least the first two 
conditions of embeddability are satisfied.  

We remark that there is another reduction of redundancy studied in \cite{bh}, which
is the  requirement that all $b$-arcs be essential: We say that a foliated disc
is {\bf essential} if a negative vertex is never attached by a saddle to a positive
vertex which is adjacent to it in the vertex string. If it was inessential, we could
reduce the foliation of the disc by eliminating that negative vertex and the saddle
attaching it to an adjacent positive vertex, without altering the embeddability of the
disc and the boundary braid. So from now on we will discard all inessential discs.
Since the test for essentiality is very easy, this is a very inexpensive way to reduce
redundancy in our list.  After we have tested the first two embeddability conditions
and eliminated inessential $b$-arcs we will run the third text for embeddability in
Theorem \ref{embeddability theorem}. 

In order to keep track of the information necessary for a possibly essential embeddable
insertion, we have invented \textbf{another combinatorial description of the disc $D$},
also suitable for our implementation.
We will code our disc $D$ by two sets: the cyclically ordered set of
boundary points of saddles, read counterclockwise along $\partial D$,
 and the set of $bb$-saddles. The reader might like to compare what follows with
Figure \ref{thedisc}.
Each $aa$-saddle has two points on the boundary $\partial
D$, that we will call \textbf{points of type Q}; each $ab$-saddle has one point
on $\partial D$ that we will call \textbf{points of type R}; and each
$bb$-saddle  lies entirely in the interior of $D$, with no point on $\partial D$. 

With each boundary point we will associate: the sign of the 
attached saddle (overline points corresponding to negative saddles); 
the ordered list of two or three vertices to which the saddle 
is attached, and the ordinal number of the saddle. 

The double index of $Q_{i,j}$ (or $\ol{Q}_{i,j}$)
is such that the cyclic counterclockwise order about the saddle 
of the boundary point and the two vertices is $i,Q,j$.

The triple index of $R_{i,v,j}$ (or $\ol{R}_{i,v,j}$) is such
that the cyclic clockwise order of this point and the three vertices around the
saddle is $R,i,v,j$ (in particular, the central index is the negative vertex).

We call \textbf{initial and final vertex} of a boundary point respectively
the first and the last positive vertex, as they occur as indices of the point. In the
sequence read along the boundary, two consecutive points always have
the final vertex of the preceding point equal to the initial vertex of the following point.

Each $bb$-tile can be coded by: the sign of the saddle, the
ordered sequence of the four vertices around it, and its ordinal number.

These combinatorial data are clearly in bijective correspondence
with the set of data given in Section \ref{drawthedisc}, so they are
sufficient to draw the foliated disc and to 
read the extended boundary word (cf \cite{bh}). 

For instance, the disc of Section \ref{drawthedisc} is described by (cf Figure \ref{thedisc}) 
$$\partial D=\{Q_{1,5}(11),Q_{5,6}(10),R_{6,4.1,7}(1),\ol{R}_{7,0.1,3}(6),R_{3,0.1,8}(13),Q_{8,4}(5),$$
$$\ol{R}_{4,8.1,11}(14),\ol{R}_{11,8.2,10}(15),Q_{10,9}(2),Q_{9,10}(2),R_{10,8.2,11}(3),R_{11,8.1,4}(4),$$
$$Q_{4,8}(5),\ol{R}_{8,0.1,7}(7),\ol{R}_{7,0.2,6}(9),Q_{6,5}(10),Q_{5,1}(11)\};$$
$$bb\mbox{-saddles}=\{\ol{[6,0.2,2,4.1]}(8),[7,0.2,2,4.1](12)\}.$$

\ 

We code an insertion arc $\alpha$ by listing all the consecutive points of type 
\textbf{Q} and \textbf{R} contained in $\alpha$. The position of the new inserted negative 
vertex $v$ can be found as follows:
\bi
\item for each point $Q_{i,j}$ or $R_{i,u,j}$ in $\alpha$: if $j-i\leq2$(mod $n$), then the insertion 
is inessential; otherwise we get $i+1<v<j-1$;
\item for each point $\ol{Q}_{i,j}$ or $\ol{R}_{i,u,j}$ in $\alpha$: if $i-j\leq2$(mod $n$), 
then the insertion is inessential; otherwise we get $j+1<v<i-1$.
\ei

These conditions correspond to the cyclic order of vertices around saddles.
This means also that if $\alpha$ contains more than one point, we have to
intersect conditions coming from different points: if the intersection is
empty, the insertion is not embeddable.  We can list all possible intervals of
insertion in order of increasing length (that is the number of singular
boundary points), giving for each of them the possible essential position of
the new negative vertex. If some of the two preceding conditions eliminate some
arc $\alpha$ of length $k$, then all arcs with length greater than $k$ and
containing $\alpha$ are inessential or not embeddable for the same reason. The
complete list is done by our GAP procedure GetInsertionArcs($n,\partial D$).

In our example (see Figure \ref{thedisc}), an insertion arc is for instance 
$\alpha=\{\ol{R}_{7,0.1,3},$ $R_{3,0.1,8}\}$, with 
$4<v<6$ from the first point and $4<v<7$ from the second, hence it must be $4<v<6$.

An insertion arc $\alpha$ has an \textbf{initial} and a \textbf{final vertex} 
$(i(\alpha),f(\alpha))$, respectively the initial
vertex of its first point, and the final vertex of its last point. In our example they are $(7,8)$.


\begin{figure}
\begin{center}\includegraphics{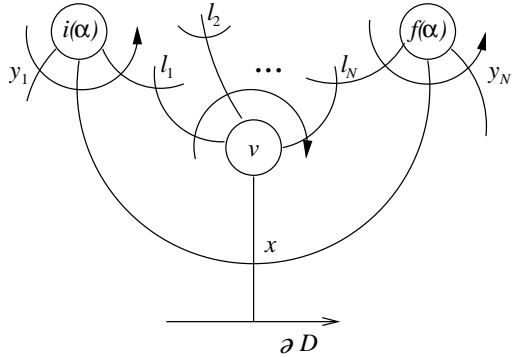}\end{center}
\caption{How to find the possible level of the new saddle.}
\label{levelofsaddle}
\end{figure}

To know at which level $x$ the new saddle can be located, we need to look at levels of 
points of $\alpha$ and to consider the cyclic order of all saddles around the new vertex,
and around the two positive vertices $i(\alpha),f(\alpha)$ (see Figure \ref{levelofsaddle}).

Suppose $\alpha=\{P_1(l_1),P_2(l_2),\ldots P_N(l_N)\}$, where each $P_j$ can be either
of \textbf{Q} or of \textbf{R} type, either positive or negative, and each $l_j$ indicates
its level. Then we must have $l_1<l_2<\cdots<l_N<x$ (in cyclic order). 
Now we have to look for the point on $\partial D$ immediately preceding $P_1$ and containing 
the index $i(\alpha)$: suppose $y_1$ is its level: then we must have $y_1<x<l_1$.
Also we have to look for the point on $\partial D$ immediately following $P_N$ and containing 
the index $f(\alpha)$: suppose $y_N$ is its level: then we must have $l_N<x<y_N$.

In our example: $6<13<x$, that means $13<x<6$, is the first requirement; $l_1=1$ and $l_N=5$,
therefore $1<x<6$ and $13<x<5$. Intersecting these cyclic intervals we get $1<x<5$.

The sign of the new saddle only depends on the relative position of $i(\alpha),v,f(\alpha)$:
\bi
\item if $f(\alpha)<v<i(\alpha)$ then $s$ is positive;
\item if $i(\alpha)<v<f(\alpha)$ then $s$ is negative;
\ei
therefore the possible range of $v$ obtained before can be divided in two parts, giving different
signs for $s$.
Our GAP procedure GetSaddles($n,P+N-1,\partial D)$ gives all these results.

In our example, $f(\alpha)=8<4<v<6<7=i(\alpha)$ hence $s$ is positive.

\ 

At this stage, before proceeding with the expensive $b$-arcs test, we can
perform an easy \textbf{permutation test}: when we do an insertion, 
it is easy to see how the corresponding 
extended word changes: a new saddle $[\pm1,[i(\alpha),v,f(\alpha))],x]$ is inserted,  
this corresponds to inserting a new letter $(i(\alpha),f(\alpha))^{\pm1}$ 
(or $(f(\alpha),i(\alpha))^{\pm1}$) at the same level of $EW$: we get a longer word, which still
must satisfy conditions given in Proposition \ref{extended word} for an extended word.
So we can try all possible combinations of $[\pm1,[i(\alpha),v,f(\alpha))],x]$ and check the 
corresponding permutation, to discard the impossible ones.

For instance our possible insertions for the chosen arc $\alpha$ are $[+1,[7,v,8],x]$, with
$x$ ranging between 1 and 5. So the new word might be one of the following:
$$EW(D_1)=(7,6)\mathbf{(8,7)}(10,9)(11,10)(11,4)(8,4)\ol{(7,3)}\ol{(8,7)}\ol{(6,2)}\cdot$$
$$\cdot\ol{(7,6)}(6,5)(5,1)(7,2)(8,3)\ol{(11,4)}\ol{(11,10)};$$
$$EW(D_2)=(7,6)(10,9)\mathbf{(8,7)}(11,10)(11,4)(8,4)\ol{(7,3)}\ol{(8,7)}\ol{(6,2)}\cdot$$
$$\cdot\ol{(7,6)}(6,5)(5,1)(7,2)(8,3)\ol{(11,4)}\ol{(11,10)};$$
$$EW(D_3)=(7,6)(10,9)(11,10)\mathbf{(8,7)}(11,4)(8,4)\ol{(7,3)}\ol{(8,7)}\ol{(6,2)}\cdot$$
$$\cdot\ol{(7,6)}(6,5)(5,1)(7,2)(8,3)\ol{(11,4)}\ol{(11,10)};$$
$$EW(D_4)=(7,6)(10,9)(11,10)(11,4)\mathbf{(8,7)}(8,4)\ol{(7,3)}\ol{(8,7)}\ol{(6,2)}\cdot$$
$$\cdot\ol{(7,6)}(6,5)(5,1)(7,2)(8,3)\ol{(11,4)}\ol{(11,10)}.$$
For them we find the following permutations:
$$\rho(EW(D_j))=(1,5,6,10,9)(2)(3)(4)(7)(8)(11)$$
(the same for all, since $(8,7)$ commutes with second, third and fourth letters),
which satisfies all conditions required by Proposition \ref{extended word}.

\ 

Now we have to perform the test for $b$-arcs.
For this we need both $\partial D$ and the set of $bb$-saddles.
Also, we first need to see \textbf{how the code changes after an insertion}.

\bi
\item if $i<v<i+1$ and in the same interval there are other existing
negative vertices, we have to decide (by the embeddability test)
the exact order of them in this interval, and rename old vertices in it
if necessary;
\item if $y-1<x<y$, put $x=y$ and for any saddle which was at level $z\geq y$
put it at level $z+1$;
\item substitute $\alpha$ by $R^{sign(s)}_{i(\alpha),v,f(\alpha)}$;
\item for each point $Q^\varepsilon_{i,j}$ in $\alpha$, its other corresponding
$Q^\varepsilon_{j,i}$ becomes $R^\varepsilon_{j,v,i}$;
\item each point $R^\varepsilon_{i,u,j}$ in $\alpha$ becomes a $bb$-saddle
$[\varepsilon,[i,u,j,v]]$.
\ei
If for instance we make the insertion of $[+1,[7,4.1,8],4]$ on the chosen $\alpha$, 
the description of the new disc is (see Figure \ref{newdisc}):
$$\partial D=\{Q_{1,5}(12),Q_{5,6}(11),R_{6,4.2,7}(1),R_{7,4.1,8}(4),Q_{8,4}(6),$$
$$\ol{R}_{4,8.1,11}(15),\ol{R}_{11,8.2,10}(16),Q_{10,9}(2),Q_{9,10}(2),R_{10,8.2,11}(3),R_{11,8.1,4}(5),$$
$$Q_{4,8}(6),\ol{R}_{8,0.1,7}(8),\ol{R}_{7,0.2,6}(10),Q_{6,5}(11),Q_{5,1}(12)\};$$
$$bb\mbox{-saddles}=\{\ol{[0.1,3,4.1,7]}(7),\ol{[6,0.2,2,4.2]}(9),$$
$$[7,0.2,2,4.1](13),[0.1,3,4.1,8](14)\}.$$

\begin{figure} 
\begin{center}\includegraphics{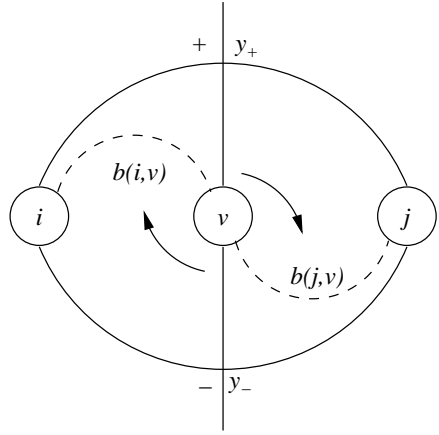}\end{center}
\caption{How to find the $b$-arcs about a negative vertex.}
\label{barcs}
\end{figure}

\noindent\textbf{How to read the $b$-arcs:} if a negative vertex $v$ (yet existing or 
newly inserted) is attached to only two saddles
(for the new vertex, this corresponds to an insertion arc of length 1), they (see Figure \ref{barcs})
will surely be one positive (say at level $y_+$) and one negative saddle (say at level $y_-$),
and they will be connected to the same two positive vertices, say $i$ and $j$.
If $i<v<j$, then in the interval $(y_+,y_-)$ we have the $b$-arc $b(j,v)$ and
in the interval $(y_-,y_+)$ we have the $b$-arc $b(i,v)$.

If a negative vertex $v$ (yet existing or newly inserted) is attached to more than two saddles
(for the new vertex, this corresponds to an insertion arc of length greater than or equal to 2), 
we have to list these saddles in cyclic order about $v$ (in our example, they are
$[7,4.1,8](4)$, $\ol{[7,0.1,3,4.1]}(7)$, $[3,0.1,8,4.1](14)$); between any two consecutive
 of them, 
find the $b$-arc between $v$ and the only other positive vertex which is in common for the two saddles
(in our example: in $(4,7)$ have $b(7,4.1)$; in $(7,14)$ have $b(3,4.1)$; in $(14,4)$ have 
$b(8,4.1)$).

As explained in \cite{bh}, we have to compose an array with: in the first column
the $P+N-1$ intervals $(k,k+1)$ between two consecutive saddles; in the following
$N$ columns, the $b$-arcs for each negative vertex in the corresponding intervals;
in the last column, the $gb$-arcs (they are as many as the $aa$-saddles).
For instance, for our new disc of Figure \ref{newdisc} we see the array of Figure 
\ref{barcstable}.

\begin{figure}
\begin{center}
$$\ba{cccccccc}
(1,2)   & b(3,0.1) & b(2,0.2) & b(8,4.1) & b(6,4.2) & b(4,8.1)  & b(11,8.2) & gb(9,10) \\
(2,3)   & b(3,0.1) & b(2,0.2) & b(8,4.1) & b(6,4.2) & b(4,8.1)  & b(11,8.2) &  \\
(3,4)   & b(3,0.1) & b(2,0.2) & b(8,4.1) & b(6,4.2) & b(4,8.1)  & b(10,8.2) &  \\
(4,5)   & b(3,0.1) & b(2,0.2) & b(7,4.1) & b(6,4.2) & b(4,8.1)  & b(10,8.2) &  \\
(5,6)   & b(3,0.1) & b(2,0.2) & b(7,4.1) & b(6,4.2) & b(11,8.1) & b(10,8.2) & gb(4,8) \\
(6,7)   & b(3,0.1) & b(2,0.2) & b(7,4.1) & b(6,4.2) & b(11,8.1) & b(10,8.2) &  \\
(7,8)   & b(7,0.1) & b(2,0.2) & b(3,4.1) & b(6,4.2) & b(11,8.1) & b(10,8.2) &  \\
(8,9)   & b(8,0.1) & b(2,0.2) & b(3,4.1) & b(6,4.2) & b(11,8.1) & b(10,8.2) &  \\
(9,10)  & b(8,0.1) & b(6,0.2) & b(3,4.1) & b(2,4.2) & b(11,8.1) & b(10,8.2) &  \\
(10,11) & b(8,0.1) & b(7,0.2) & b(3,4.1) & b(2,4.2) & b(11,8.1) & b(10,8.2) & gb(5,6) \\
(11,12) & b(8,0.1) & b(7,0.2) & b(3,4.1) & b(2,4.2) & b(11,8.1) & b(10,8.2) & gb(1,5)\\
(12,13) & b(8,0.1) & b(7,0.2) & b(3,4.1) & b(2,4.2) & b(11,8.1) & b(10,8.2) &  \\
(13,14) & b(8,0.1) & b(2,0.2) & b(3,4.1) & b(7,4.2) & b(11,8.1) & b(10,8.2) &  \\
(14,15) & b(3,0.1) & b(2,0.2) & b(8,4.1) & b(7,4.2) & b(11,8.1) & b(10,8.2) &  \\
(15,16) & b(3,0.1) & b(2,0.2) & b(8,4.1) & b(7,4.2) & b(4,8.1)  & b(10,8.2) &  \\
(16,1)  & b(3,0.1) & b(2,0.2) & b(8,4.1) & b(7,4.2) & b(4,8.1)  & b(11,8.2) &  \\
\ea$$
\end{center}
\caption{Table of $b$- and $gb$-arcs.}
\label{barcstable}
\end{figure}

\

\noindent\textbf{Remark:} When an $aa$-saddle occurs, no $b$-arc is changed; when an
$ab$-saddle occurs, only changes the $b$-arc connected to the negative vertex involved
in the
$ab$-saddle; when a $bb$-saddle occurs, only change the two $b$-arcs connected to the
negative vertices  involved in the $bb$-saddle.

\

It is necessary that along each row of the array the arcs which appear do not
interlock with each other. All these checks are performed by our GAP procedure
Embeddable$(n,P+N-1,\partial D,bb$-saddles). 
Finally, our GAP procedure InsertVertices$(n,P+N-1,\partial D,bb$-saddles,
depth) recursively tries all possible insertions of $ab$-tiles and returns all
the resulting new embeddable essential tiled disc.

\section{How to get the boundary word.\\Changes in foliation.}
\label{boundaryword}
To get the boundary word $BW$ we have to eliminate from $EW$ the
$N$ strands corresponding to the $N$ unlinked circles
about each negative vertex. In this way, from each letter of $EW$ we
will get a new letter or subword of the boundary braid $BW$, which is
a braid of $B_{P-N}$.

To do so, let us introduce some useful braids (called \emph{descending cycles} 
in \cite{bkl:1}):
\bd When $p>q$, call $\delta_{p,q}$ the braid $(p,p-1)(p-1,p-2)\cdots(q+1,q)$.\ed
Notice that $\delta_{p,p-1}=(p,p-1)=\sigma_{p-1}$. The associated permutation is
$\rho(\delta_{p,q})=(q,q+1,\ldots,p)$.

\


\begin{figure}
\begin{center}\includegraphics{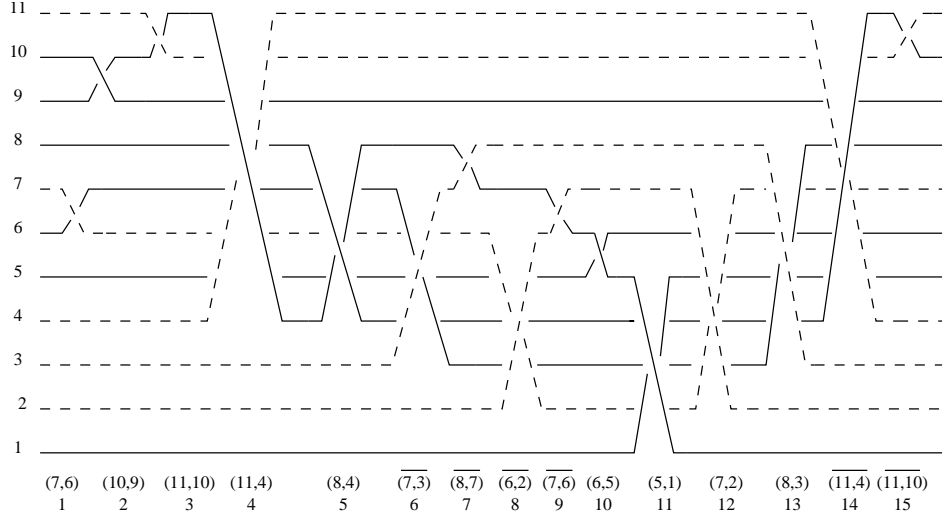}\end{center}
\caption{The braid $EW$. The dashed strands must 
be removed to find the boundary braid.}
\label{extendedbraid}
\end{figure}

\begin{figure}
\begin{center}\includegraphics{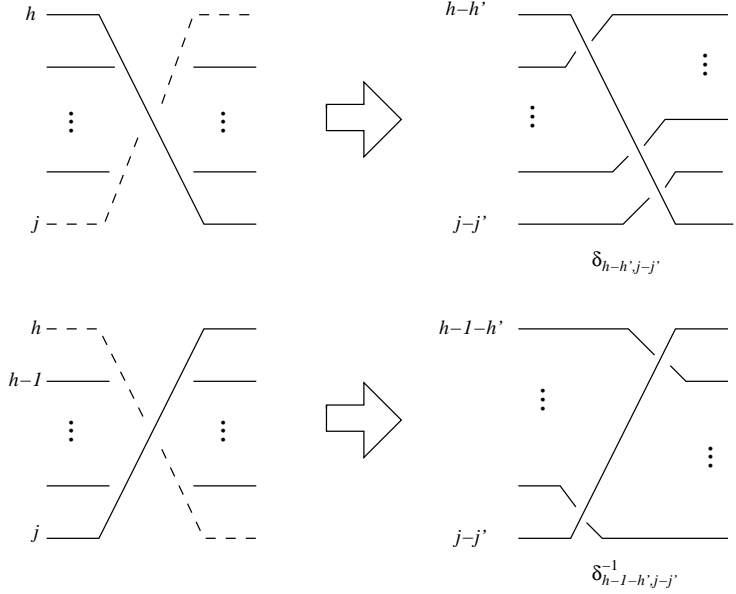}\end{center}
\caption{The braids $\delta$'s obtained from a 
generator by removing one strand.}
\label{deltas}
\end{figure}

The strands to be eliminated are numbered, at the beginning of $EW$, by the 
$N$ 1-cycles of $\rho(EW)$.
At each letter of $EW$ their numbers might change, as can be seen for instance 
in Figure \ref{extendedbraid}. The rule is: suppose $L_k=\{i_1,i_2,\ldots i_N\}_k$ is the list of
levels of strands to be eliminated just before the $k^{th}$ letter of $EW$, and this
letter is $(h,j)^\varepsilon$ (see Figure \ref{deltas}):
\be
\item if both $h,j$ are not in $L_k$, and $h',j'$ are the numbers
of elements of $L_k$ which are less than $h,j$ respectively; then get 
$(h-h',j-j')^\varepsilon$; $L_{k+1}=L_k$;
\item if $j\in L_k$ and $h$ is not in $L_k$: if $h>j+1$, then get $\delta_{h-h',j-j'}$;
if $h=j+1$, then get the empty word $e$; $L_{k+1}=(L_k\setminus\{j\})\cup\{h\}$;
\item if $h\in L_k$ and $j$ is not in $L_k$: if $h>j+1$, then get $\delta^{-1}_{h-1-h',j-j'}$;
if $h=j+1$, then get $e$; $L_{k+1}=(L_k\setminus\{h\})\cup\{j\}$;
\item if both $h,j\in L_k$, then get $e$; $L_{k+1}=L_k$.
\ee

In our example we get 
$$BW=e(6,5)e\delta_{6,2}(5,2)\delta_{4,2}eee(5,4)(4,1)e\delta^{-1}_{5,2}\delta^{-1}_{6,2}e=$$
$$=(6,5)\delta_{6,2}(5,2)\delta_{4,2}(5,4)(4,1)\delta^{-1}_{5,2}\delta^{-1}_{6,2}.$$

Our GAP procedure BoundaryBraid$(EW,P,N)$ computes the boundary braid of a disc with
extended word $EW$ and $(P,N)$ vertices.

\ 


In what follows we explain how the defining relations and inversion in the braid group, 
performed on the extended word, affect the topology of the foliated disc and the corresponding 
boundary word. 
\begin{question}
How do these changes in foliation affect the possibility of inserting $ab$-tiles?
\end{question}
To describe what happens let us first notice some properties of the $\delta_{i,j}$'s.

\bp 
\be 
\item $\delta_{i,j}\delta_{h,k}=\delta_{h,k}\delta_{i,j}$ if $i>j>h>k$;
\item $\delta_{i,j}\delta_{h,k}=\delta_{h-1,k-1}\delta_{i,j}$ if $i>h>k>j$;
\item $(i,j)\delta_{h,k}=\delta_{h,k}(i+1,j+1)$ if $h>i>j>k$;
\ee
\ep

In what follows, when we say `letter' we mean letters in band generators or $\delta$'s.

\ 

The easiest `relation' we want to describe is free reduction $(i,k)\ol{(i,k)}$.

\bp Let $EW$ be the extended word of an embeddable disc $D$ with $(P,N)$
vertices. If $EW$ is reducible (i.e. it has two
consecutive letters, one of which is the inverse of the other), then the
reduced word $EW'$ is the extended word of another embeddable disc $D'$ 
with $(P-1,N-1)$ vertices. Moreover, the two boundary braids $BW,BW'$ define the same
braid.
\ep

\begin{figure}
\begin{center}\includegraphics{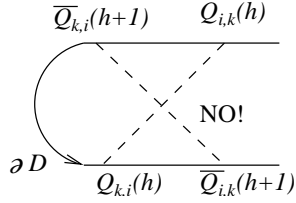}\end{center}
\caption{Two impossible $aa$-saddles.}
\label{impossibleaasaddles}
\end{figure}

\begin{figure}
\begin{center}\includegraphics{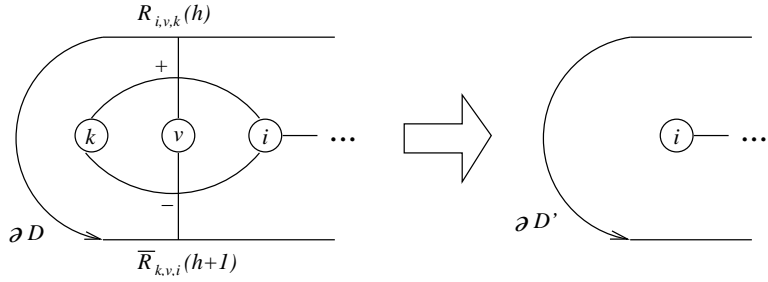}\end{center}
\caption{Elimination of two adjacent inverse saddles.}
\label{elimination}
\end{figure}

\begin{figure}
\begin{center}\includegraphics{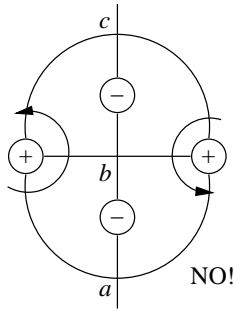}\end{center}
\caption{An impossible $bb$-saddle.}
\label{impossiblebbsaddle}
\end{figure}

\noindent\textbf{Proof:} Suppose that $w_hw_{h+1}=(i,k)\ol{(i,k)}$. Then, along 
$\partial D$ we see either the two pairs of
consecutive points $Q_{i,k}(h)\ol{Q}_{k,i}(h+1)$ and
$Q_{k,i}(h)\ol{Q}_{i,k}(h+1)$, or one of the pairs
$R_{i,v,k}(h)\ol{R}_{k,v,i}(h+1)$ or $R_{k,v,i}(h)\ol{R}_{i,v,k}(h+1)$.
The first case is impossible (see Figure \ref{impossibleaasaddles}), 
since the two saddles would cross each
other. In the second case, the positive vertex $k$ (resp. $i$) is
isolated from any other $aa$- or $ab$-saddle (see Figure \ref{elimination}); 
there cannot be other saddles involving $k$, because the two points on the boundary
are consecutive, and no other $bb$-saddle can occur among the two vertices, 
to respect embeddability (see Figure \ref{impossiblebbsaddle}). 
So we can reduce $EW$ by deleting the two inverse letters,
to get a braid word in $B_{P-1}$ representing a disc
$D'$ (see Figure \ref{elimination}) which differs from $D$ by not having the two
corresponding saddles, the vertices $k$ (resp. $i$) and $v$.
The eliminated positive vertex, because of its position in $D$,
is one of those in the one-cycles of the permutation. 
The new word $EW'$ has length
$(P-1)+(N-1)-1$ and one less one-cycle, so it still has a good permutation. 
We have eliminated two $ab$-saddles from $D$, so we have
eliminated some $b$-arcs, therefore since $D$ was embeddable,
so is $D'$.

The two boundary braids live in the same braid group, since $P-N=P-1-(N-1)$.  
Moreover the two consecutive inverse letters of $EW$ go into two
consecutive inverse letters of $BW$ (maybe $e$), so $BW=BW'$ as braid elements.
\cvd

\


\noindent\textbf{Remark:} The inverse operation, ie the insertion in an extended word of
a pair of inverse letter, is not always admissible: first of all, one of the two indices
of the inserted letters must be new, or we have to add one to all indices after this.
Even so, the insertion might produce an inessential or non embeddable disc: for instance
given $EW=(6,5)(5,4)\ol{(4,2)}(3,1)(5,3)$, if we insert the pair $(7,2)\ol{(7,2)}$ between
the second and the third letter, we get an inessential disc. If we insert in the same
position the pair $(7,3)\ol{(7,3)}$, we get an  impossible disc.

\bp Let $EW$ be the extended word of an embeddable disc $D$. 
If in  $EW$ two consecutive letters commute, then the
related word $EW'$ is the extended word of another embeddable disc $D'$ 
such that the two boundary braids $BW,BW'$ are the same braid.
\ep

\begin{figure}
\begin{center}\includegraphics{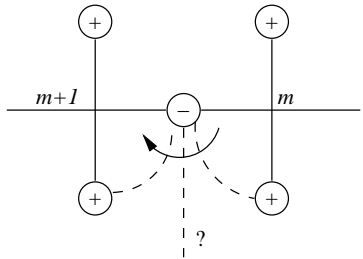}\end{center}
\caption{Two impossible $ab$-saddles.}
\label{impossibleabsaddles}
\end{figure}

\noindent\textbf{Proof:}  Suppose that $w_mw_{m+1}=(i,j)^{\varepsilon}(h,k)^{\eta}$, with
$\varepsilon,\eta\in\{\pm1\}$. Then the two
corresponding  saddles cannot have any common vertex. For, if there is a common vertex 
it is negative. But then it is impossible to foliate regularly the
region about this negative vertex which is bounded by the two 
 saddles, see Figure \ref{impossibleabsaddles}.

Therefore, exchanging the order of the two saddles does not affect
the cyclic order of saddles about vertices, nor of vertices about saddles.
In the interval $(m-1,m+2)$ only these two saddles occur, so only
their $b$- or $gb$-arcs change their occurrence: but since they
do not share any common vertex, and all the other arcs remain unchanged,
this commutation does not affect the $b$-arc test, so $D'$ is still
embeddable.

Clearly the only change in $BW$ is the commutation of the two
corresponding letters. This commutation leads to a different braid word
representing the same braid.
\cvd

\bp\label{inversion} Let $EW$ be the extended word of an essential embeddable disc $D$. 
Then $EW^{-1}$ is the extended boundary word of an essential embeddable
disc $\ol{D}$, which can be drawn as follows: look at $D$ from its negative
side, let all names of the $(P,N)$ vertices unchanged, change sign to all
saddles, and change numbers of saddles by the reversing permutation
$(1,P+N-1)(2,P+N-2)\ldots$
\ep

\noindent\textbf{Proof:} For this proof, we use results of section \ref{section:The H-theta
sequence}. $D$ is essential and embeddable if and only if it has an 
essential and embeddable $H_\theta$-sequence. Consider the sequence read from 
it in the reverse order: so each $ab$- or $bb$-saddle changes its sign; assign 
opposite sign also to each $aa$-saddle. The resulting sequence is clearly essential
and embeddable as the previous one, and it corresponds to the inverse
extended boundary word: it has inverse induced permutation, hence a good
one again. If we draw the disc starting from $EW^{-1}$ we get saddles
in the reverse order. \cvd

\

An example is shown in Figure \ref{discogretchen}.

\bp Let $EW$ be the extended word of an embeddable disc $D$. 
If in  $EW$ we can perform a relation between two consecutive letters
sharing an index, then the
related word $EW'$ is the extended word of another embeddable disc $D'$ 
such that the two boundary braids $BW,BW'$ are the same braid.
\ep

\begin{figure}
\begin{center}\includegraphics{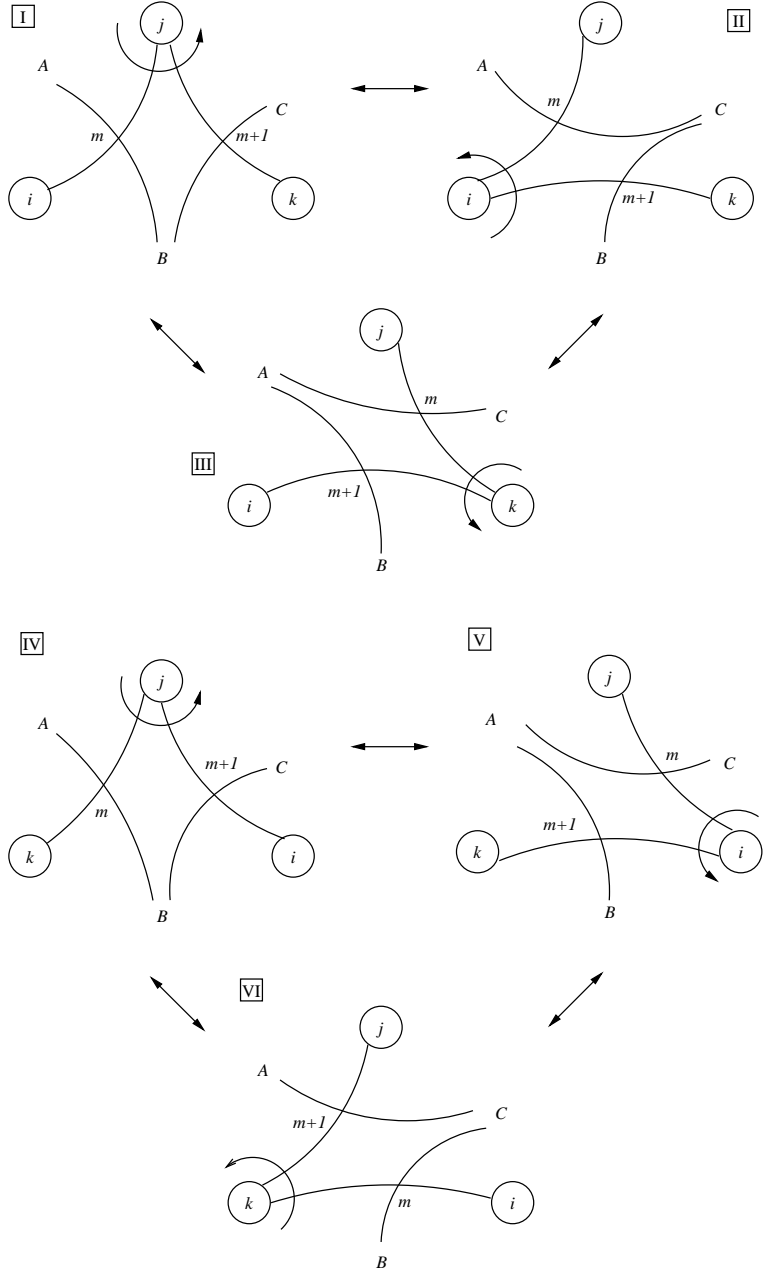}\end{center}
\caption{Change in foliation: the possible signs of saddles are: 
$I(+,+)\leftrightarrow II(+,+)\leftrightarrow III(+,+)$; 
$I(+,-)\leftrightarrow II(-,+)$; $II(+,-)\leftrightarrow III(-,+)$;
$III(+,-)\leftrightarrow I(-,+);$ 
$IV(-,-)\leftrightarrow V(-,-)\leftrightarrow VI(-,-)$; 
$IV(+,-)\leftrightarrow V(-,+)$; $V(+,-)\leftrightarrow VI(-,+)$;
$VI(+,-)\leftrightarrow IV(-,+)$.}
\label{relations}
\end{figure}

\begin{figure}
\begin{center}\includegraphics{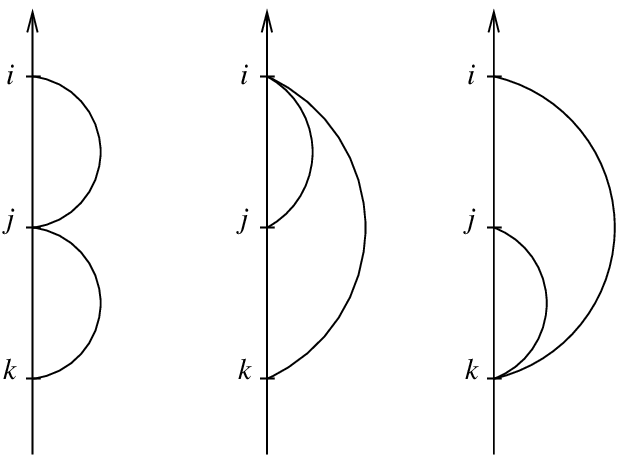}\end{center}
\caption{Relations that do not change the possibility of drawing
$b$-arcs.}
\label{connectedcomponents}
\end{figure}

\noindent\textbf{Proof:} The configurations on $D,D'$ of the possible 
relations are given in Figure \ref{relations}. 
If all three saddles are of $aa$-type, so that none of the points $A,B,C$ is a negative
vertex, then the embeddability is unchanged: in fact, in the interval 
$(m-1,m+2)$ only these two $aa$-saddles occur, and the two corresponding 
$gb$-arcs change in such a way that the possibility of drawing
$b$-arcs in their complement does not change (see Figure
\ref{connectedcomponents}).

If some of the points $A,B,C$ are negative vertices, we also have to 
take account of the change of $b$-arcs. In this case too, by examining all possible
cases, we see that the possibility
of drawing the other $b$-arcs is not affected by the relation.

The change in the boundary word is a corresponding relation,
possibly between $\delta$'s or between some $\delta$ and
some band generator, but always giving a related boundary
word representing the same braid element.
\cvd

\ 

\noindent\textbf{Remark:} Relations do not affect embeddability,
but they might affect essentiality. Moreover,
 it is true that they change the topology of the embeddable disc,
so that, for instance, they can affect the possibility of inserting an $ab$-tile.
In fact, they always change the possible insertion arcs.
This leads to the following:

\begin{question} We phrase this question as a conjecture: The list of good words to be
used to find new words by insertions of $ab$-tiles cannot be reduced by performing
relations of the braid group. \end{question}

Here is an example:
$$W_1=(65)(54)\ol{(42)}(31)(53); W_2=(65)\ol{(52)}(54)(31)(53);$$
these two words represent the same braid, since they only differ by
a relation. They are both good words, therefore they represent
embeddable positive discs. But when we examine them, we find that
in $D_1$ we can perform two embeddable essential insertions:
\be
\item along $\ol Q_{2,4}(3)$ with $6<v<1$ and positive saddle at level $4<x<5$;
\item along $Q_{5,3}(5)$ with $6<v<1$ and negative saddle at level $2<x<3$.
\ee
The extended and boundary words of the discs we get from these insertions are
respectively:
\be
\item $E_1=(65)(54)\ol{(42)}(31)(42)(53)$ and $B_1=(54)(43)(32)(31)\ol{(32)}(42)$;
\item $E_2=(65)(54)\ol{(53)}\ol{(42)}(31)(53)$ and $B_2=(54)(43)(43)\ol{(42)}(31)\ol{(43)}$.
\ee
Notice that these two 5-braids are not conjugate, since they have
different exponent sum.

\ 

But if we look at $W_2$, we find that only one essential embeddable
insertion is possible, namely $\{\ol Q_{25}(2),Q_{54}(3)\}$, 
with $6<v<1$ and positive saddle at level $4<x<5$,
from which we get the extended and the boundary word
$$E_3=(65)\ol{(52)}(54)(31)(42)(53), \, B_3=(54)(43)(32)(31)\ol{(32)}(42),$$
which is the same we got from the first insertion on $W_1$. In fact this
insertion does correpond to the other one, since the change in
foliation given by the relation has substituted the arc containing
$\ol Q_{24}$ with the arc containing now $\ol Q_{25},Q_{54}$.

\

\noindent\textbf{Remark:} We might find the other 5-braid from a different good word,
related to $W_1$ in some other way.  Or if we list our discs
using the method of $H_\theta$-sequences (see section \ref{section:The H-theta
sequence}) instead of insertions  of $ab$-tiles, we have not to be worried by these
changes in foliation.

\ 

\noindent\textbf{Change in the code after a relation (for good words):}

Suppose $W$ is a good word and we perform a relation between two
consecutive letters $w_hw_{h+1}$. Then the code for
$W$ changes as follows:

Commutation of non interlocking pairs:
$$(ij)^{\varepsilon}(kl)^{\eta}\longleftrightarrow (kl)^{\eta}(ij)^{\varepsilon}$$
(here $\varepsilon,\eta\in\{\pm1\}$).

Change in the code: exchange levels $h\longleftrightarrow h+1$ of the
four $Q$ points.

Relation between two positive letters with three ordered indices
$n\geq i>j>k\geq1$:
$$(ij)(jk)\longleftrightarrow(jk)(ik)\longleftrightarrow(ik)(ij)\longleftrightarrow$$
Change in the code: three arcs of the boundary change as follows:
$$ Q_{ij}(h)Q_{jk}(h+1)\longleftrightarrow
Q_{ik}(h+1)\longleftrightarrow Q_{ik}(h)\longleftrightarrow;$$
$$ Q_{ji}(h)\longleftrightarrow
Q_{jk}(h)Q_{ki}(h+1)\longleftrightarrow Q_{ji}(h+1)\longleftrightarrow;$$
$$ Q_{kj}(h+1)\longleftrightarrow
Q_{kj}(h)\longleftrightarrow Q_{ki}(h)Q_{ij}(h+1)\longleftrightarrow.$$

Relation between two negative letters with three ordered indices
$n\geq i>j>k\geq1$:
$$\ol{(jk)(ij)}\longleftrightarrow\ol{(ik)(jk)}\longleftrightarrow\ol{(ij)(ik)}\longleftrightarrow$$
Change in the code: three arcs of the boundary change as follows:
$$ \ol{Q}_{jk}(h)\longleftrightarrow
\ol{Q}_{jk}(h+1)\longleftrightarrow
\ol{Q}_{ji}(h)\ol{Q}_{ik}(h+1)\longleftrightarrow;$$
$$ \ol{Q}_{kj}(h)\ol{Q}_{ji}(h+1)\longleftrightarrow
\ol{Q}_{ki}(h)\longleftrightarrow \ol{Q}_{ki}(h+1)\longleftrightarrow;$$
$$ \ol{Q}_{ji}(h+1)\longleftrightarrow
\ol{Q}_{jk}(h+1)\longleftrightarrow
\ol{Q}_{ji}(h)\ol{Q}_{ik}(h+1)\longleftrightarrow.$$

Relation between one positive and one negative letter with three
cyclically ordered indices $i>j>k$:
$$(jk)\ol{(ik)}\longleftrightarrow \ol{(ij)}(jk)$$

Change in the code: three arcs of the boundary change as follows:
$$Q_{jk}(h)\ol{Q}_{ki}(h+1)\longleftrightarrow 
\ol{Q}_{ji}(h);$$
$$Q_{kj}(h)\longleftrightarrow  Q_{kj}(h+1);$$
$$\ol{Q}_{ik}(h+1)\longleftrightarrow  \ol{Q}_{ji}(h)Q_{ik}(h+1).$$

The other relation between one positive and one negative letter with three
cyclically ordered indices $i>j>k$:
$$(jk)\ol{(ij)}\longleftrightarrow \ol{(ij)}(ik)$$
Change in the code: three arcs of the boundary change as follows:
$$Q_{jk}(h)\longleftrightarrow \ol{Q}_{ji}(h)Q_{ik}(h+1);$$
$$Q_{kj}(h)\ol{Q}_{ji}(h+1)\longleftrightarrow Q_{ki}(h+1);$$
$$\ol{Q}_{ij}(h+1)\longleftrightarrow \ol{Q}_{ij}(h).$$

A similar description could be done for any extended word.

\ 

\noindent\textbf{Non conjugated good words giving the same disc after an
insertion:}

These words can be found from an essential embeddable disc with
one negative vertex, by stabilizing along one $ab$-tile (stabilizing is the
opposite process than inserting an $ab$-tile).

For instance, if the negative vertex is connected just to two
$ab$-saddles `facing each other',  of opposite sign and at
different levels (eg saddles 4 and 14 in Figure \ref{thedisc}),  
by stabilizing along one or the other we get words
with different exponent sum, hence non conjugate.

More precisely: if two words only differ for one letter $(ki)$ which is
positive and at level $s$ in the first word, but negative and at level
$t$ in the second, then if we insert an $ab$-tile with negative vertex
at level $j$, but: on $Q_{ki}(s)$ with negative saddle at level $t$ on
the first disc, and on $\ol{Q}_{ik}(s)$ with positive saddle at level
$s$ on the second, we get the same disc.

\section{The $H_\theta$-sequence}
\label{section:The H-theta sequence}
As we have seen, the test for $b$-arcs in the embeddability test is very expensive.
But there is a different method for testing embeddability, used repeatedly in the
papers of the first author and Menasco, e.g. see the proof of Lemma 3 in
\cite{bm4}.  Using it, we can avoid the expensive part of the embeddability test.

\begin{figure}
\begin{center}\includegraphics{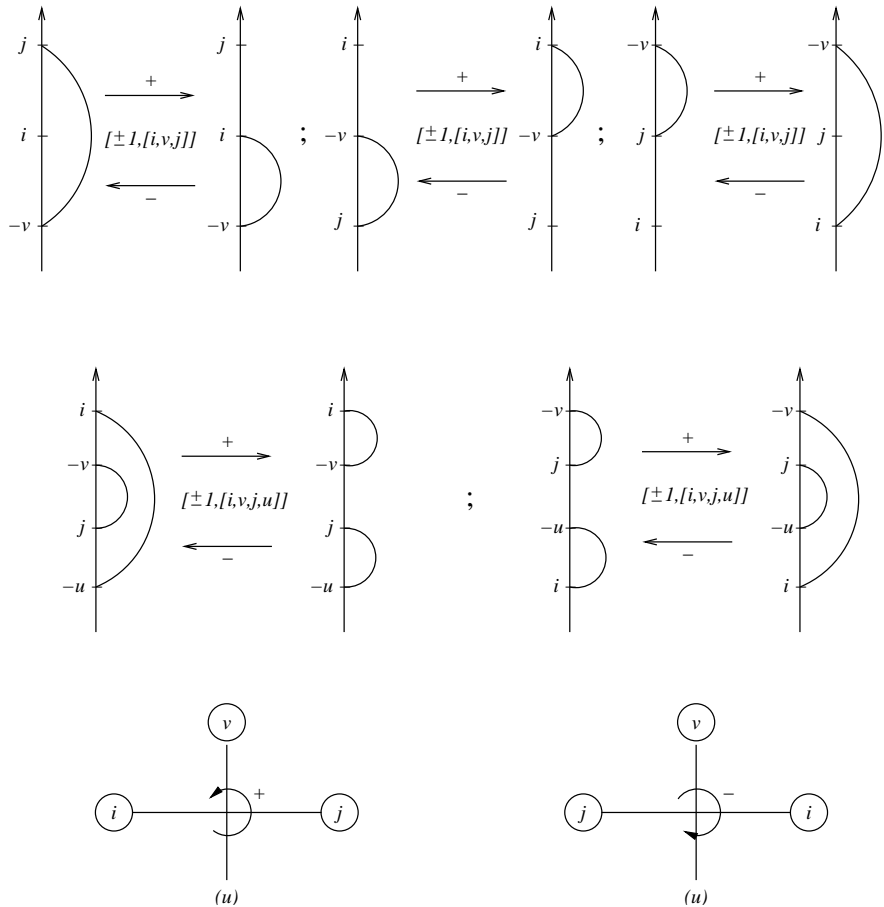}\end{center}
\caption{The $ab$-saddles and $bb$-saddles seen on the halfplanes.}
\label{signofsaddles}
\end{figure}

Look at the situation from another point of view: not at the foliation on the
disc, but at the foliations on the half-planes $H_\theta$'s running about the braid
axis $A$. In fact, almost all 
$H_\theta$'s have $P-N$ $a$-arcs and $N$ $b$-arcs (all regular leaves for $D$),
except $P+N-1$ of them, in each of which a saddle occurs. Among two consecutive
singular half-planes, all infinite regular ones are uniquely identified 
up to isotopy by the $N$ $b$-arcs. Also, the passage through one singular
half-plane is such that:
\bi
\item if an $aa$-saddle occurs, then no $b$-arc changes;
\item if an $ab$-saddle occurs: then one $b$-arc changes its positive vertex;
\item if a $bb$-saddle occurs: then two $b$-arcs exchange their positive vertices.
\ei
Moreover, the way in which these changes occur uniquely specifies the type, the names
and the sign of the $ab$- or $bb$-saddle. Only the $aa$-saddles remain unspecified
in names and sign. In Figure \ref{signofsaddles} we show all these changes and the 
corresponding saddles.

\bp A half-plane $H_\theta$ with specified vertex string $V$ with $(P,N)$
vertices is embeddable and essential if and only if the $N$ $b$-arcs
do not intersect each other, and 
each of them does not connect two consecutive vertices.
\cvd\ep

Given a vertex string, it is very easy to list all possible essential and embeddable
regular half-planes with these vertices.


\bp[cf \cite{bm4}]\label{propfilm} There is a bijective correspondence
between essential embeddable foliated discs with $(P,N)$ vertices
and $H_\theta$-sequences of length $P+N-1$ satisfying the following properties:
\bi
\item Among two consecutive half-planes in the sequence the unique change is
one corresponding to a saddle;
\item The permutation associated to the resulting cycle of saddles satisfies the 
properties of Proposition \ref{extended word};
\item No saddle occurs twice; the unique saddles involving the same pair
of positive vertices can be two $ab$-saddles of opposite sign;
\item All vertices occur at least in one saddle.	\cvd
\ei
\ep

So this is \textbf{our algorithm}: 
\be
\item For each $(P,N)$, with $P>N+1$, list all possible vertex strings, up to cyclic 
order along the axis.
\item For each vertex string, list all possible embeddable essential regular half-planes;
\item Given the complete set of half-planes found for a given vertex string, list all
the cyclic sequences of $P+N-1$ half-planes  
that satisfy all conditions of Proposition \ref{propfilm}.
\item For each found cycle, choose all possible signs of $aa$-saddles: get all possible
embeddable essential discs with $(P,N)$ vertices up to easy conjugations.
\ee

The first part of this algorithm has been implemented using the algorithm of
\cite{rs}: it corresponds to listing all \emph{necklaces} of \emph{length}
$l=P+N$, with number of \emph{colors} $k=2$ (positive and negative vertices),
and \emph{density} $d=P$ (the number of non-zero colors).  We have implemented
the algorithm as a GAP procedure, called EnumerateNecklaces($l,k,d$), which
gives us all vertex strings $V$ with $(P,N)$ vertices up to cycling. 

\begin{figure}
\begin{center}
\includegraphics{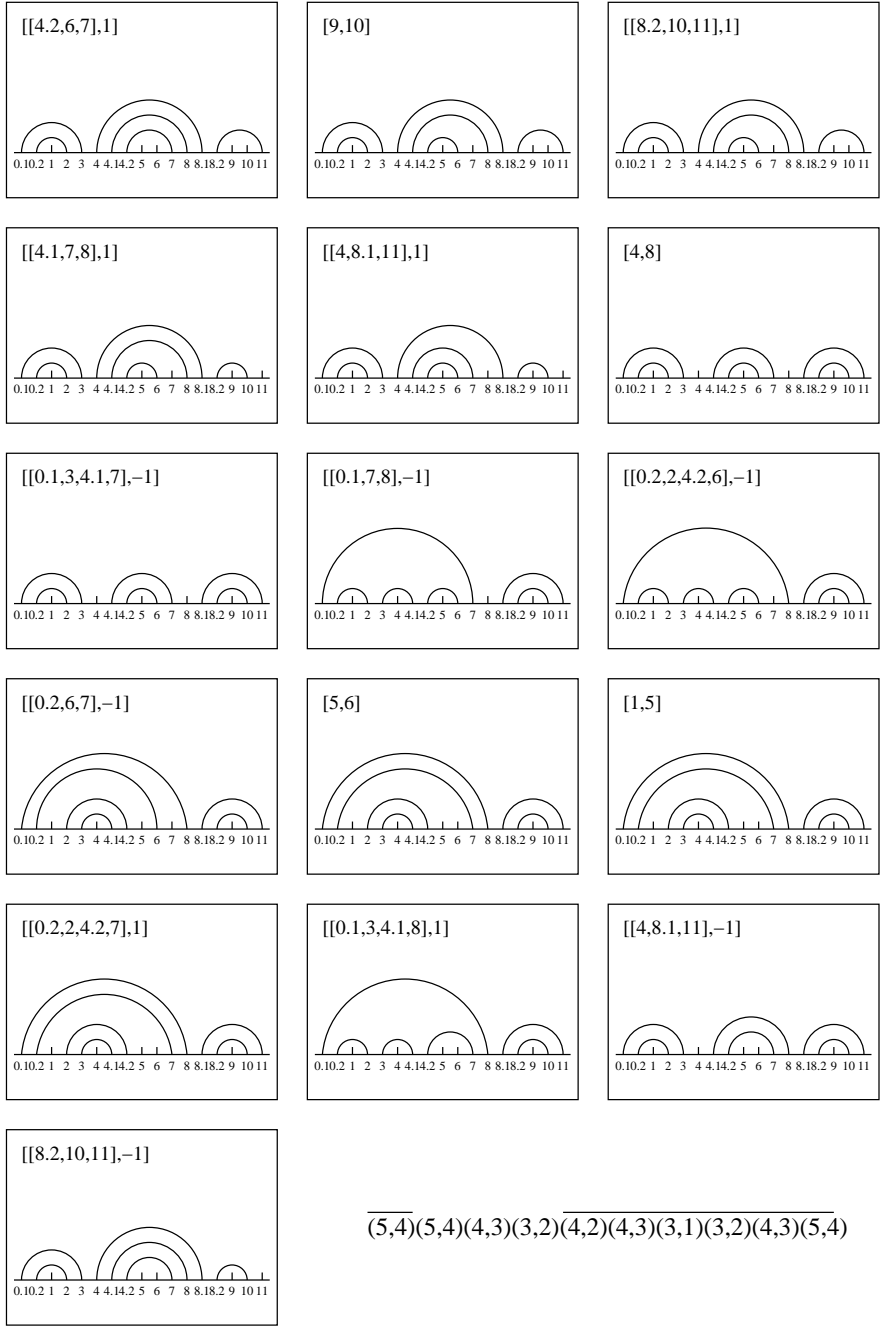}
\end{center}
\caption{The $H_\theta$-sequence of the disc of Figure \ref{newdisc}.}
\label{film}
\end{figure}


Then, the procedure EnumerateHalfPlanes($V$) lists all embeddable essential
regular half-planes in terms of their $b$-arcs. 

To list all good $H_\theta$-sequences for a given vertex string $V$, we
consider the directed graph $G$ with nodes all the regular essential embeddable
half-planes with vertex string $V$, and (directed) edges\footnote{We refrain from
using the more standard term ``arc'' for obvious reasons.} the saddles
occurring among them: notice that $aa$-saddles are loops: edges going
from one node to itself. This is made by our GAP procedure MakeGraph($V$).

Now, a good $H_\theta$-sequence corresponds to a cycle in $G$, of length $P+N-1$, 
such that the above conditions are satisfied. Such a cycle never passes twice
through the same edge, because it would pass twice through the same saddle with
the same sign. This observation is the key to our enumeration algorithm for
$H_\theta$-sequences, which is invoked by the GAP function
EnumerateCycles($G,n$), which enumerates the first $n$
$H_\theta$-sequences on the graph $G$, or all of them if $n=0$ (see the
Appendix for the details of the algorithm).




The whole process is performed in one step by our GAP procedure
ComputeCycles($P,N$), which computes all 
different
cycles starting from a given number of positive and negative vertices.


\begin{question} We noticed that many vertex strings have no associated 
$H_\theta$-sequence. This may happen for various reasons: There may be too few half-planes 
and saddles, or the resulting graph might be too disconnected, or have too few different saddles. 
 We suggest that this matter be investigated, with the goal of discarding
some vertex strings a priori, perhaps reducing considerably the running
time of the algorithm.  For example, vertex strings in which positive and negative
vertices alternate too closely have very few or no essential half-planes associated to them. 
To give an example: There are 43 vertex strings with $(8,4)$ vertices, but only 14 of
them have some cycle associated. 
\end{question}


\

\noindent\textbf{Remark:} An {\bf end-tile} on an embeddable foliated disc is a 
part of the disc which contains a positive vertex which is only
attached to one $aa$-saddle. For example, 
our disc of Figure \ref{thedisc} has two end-tiles: saddle 11 and saddle 2. These saddles
can easily be eliminated by a Markov move, also reducing the braid index. This can be
easily seen  on the extended word: if one index only occurs in one letter of $EW$, it
is surely corresponding to one end-tile; we can remove that letter and dropping all
following indices by 1, getting another word which is the extended word of a simpler
essential embeddable  disc. Our ComputeCycles($P,N$) tells us which of these 
$H_\theta$-sequences are associated to discs without end-tiles. 

\

\begin{figure}
\begin{center}
\includegraphics{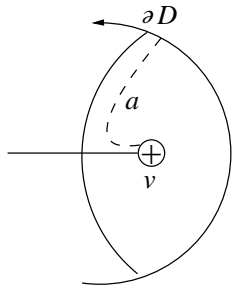}
\end{center}
\caption{A vertex of type $(a)$.}
\label{typeavertex}
\end{figure}

How does the removal of an end-tile change the corresponding boundary word? 
If it was $EW\in B_P, BW\in B_n, n=P-N$, taking away the end-tile and the corresponding strand, 
we get a new $EW'\in B_{P-1},BW'\in B_{n-1}, n-1=P-1-N$. The subword of $BW$ corresponding
to the letter of $EW$ we have deleted does not occur any more in $BW'$; the indices 
which have been dropped by 1 in $EW'$ must be dropped by 1 in $BW'$ as well.
For instance, in our example of Figure \ref{extendedbraid}, if we cancel from $EW$ the second letter
we get the new boundary word
$$BW'=\delta_{5,2}(5,2)\delta_{4,2}(5,4)(4,1)\delta^{-1}_{5,2}\delta^{-1}_{5,2}.$$

\


\noindent\textbf{How to draw the $H_\theta$-sequence:} 
We have written a GAP procedure PSFilm($C,G$), which, given a cycle $C$ for the
graph $G$, outputs a PostScript file with the drawing of the $P+N-1$
half-planes, each with its vertices and $b$-arcs, the saddle occurring between
two consecutive of them (with no sign for the $aa$-saddles) and a boundary word
resulting from arbitrarily assigning signs to $aa$-saddles. For instance, the
$H_\theta$-sequence of our disc of Figure \ref{newdisc} can be seen in Figure
\ref{film}.

There is also the possibility of exporting in GML (Graph Modelling Language) a
drawing of the tiled disc corresponding to a cycle $C$ \textit{via} the
procedure DrawDisc($C,G$). The resulting file is readable by GML-aware
software, such as Graphlet~\cite{gr}. Usually standard planar-graph layout
algorithms are able to display such tiled discs correctly, but sometimes a bit
of tweaking is required. We plan to implement in the future a more
sophisticated layout algorithm that uses the known cycling order of the
vertices.

\section{Some interesting data}
\label{sec:data}

The main theorem in \cite{bm5} is the basis for the algorithm
in \cite{bh}. This theorem was later re-proved as Theorem 4.3 of \cite{bf}, which we
now present.  Before so-doing we need several definitions.

To each tiled disc we can associate the \emph{graph of singular leaves},
in which we consider as vertices of the graph only the $(P,N)$ 
intersection points with the axis.
 
The \emph{valence} of a vertex in the graph is the number 
of singular leaves which meet at that vertex. Each non-singular leaf which 
has an endpoint at the vertex is necessarily type $a$ or type $b$, with 
the type of that leaf changing only after the passage through a singular leaf. 
We define the \emph{type} of the
vertex to be the cyclic array of $a$'s and $b$'s which describes the
non-singular leaf types as we travel around the vertex in the order in
which they are encountered in the fibration. 
The \emph{sign} of a vertex (as vertex of the graph) is the
cyclic array of signs of the singular leaves as we travel around the vertex,
again ordered by the order in which they are encountered in the
fibration. The main theorem in \cite{bm5}, which is also Theorem 4.3 of 
\cite{bf}, asserts:

\begin{figure}
\begin{center}
\includegraphics{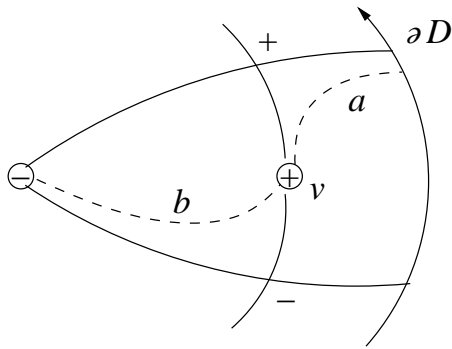}
\end{center}
\caption{A vertex of type $(a,b)$ and sign $(+,-)$.}
\label{typeabvertex}
\end{figure}

\begin{figure}
\begin{center}
\includegraphics{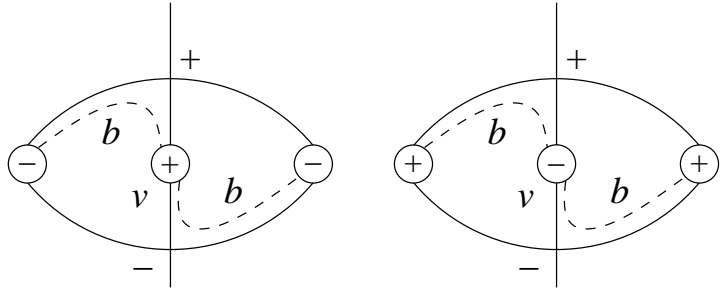}
\end{center}
\caption{Two vertices of type $(b,b)$ and sign $(+,-)$.}
\label{typebbvertex}
\end{figure}

\begin{theorem}
\label{theorem:simplifying the disc}
Let $D$ be an embeddable disc which supports a braid foliation. Then there
is a sequence of embeddable foliated discs:
$$D = D_1\to D_2\to\cdots D_k$$
such that $D_k$ has a radial foliation, without singularities, and the
graph of singular leaves for $D_{i+1}$ is obtained from that for $D_i$ 
by one of the following:
\begin{enumerate}
\item The graph of singular leaves for $D_i$ contains a vertex $v$ of
valence 1 (see Figure \ref{typeavertex}).  Delete $v$ and the 
unique singular leaf which ends at $v$.
\item The graph of singular leaves for $D_i$ contains a vertex of valence
2, type $(a,b)$ and sign $(+,-)$ (see Figure \ref{typeabvertex}). Do an
$ab$-exchange move, as defined in \cite{bf}. This move deletes two
vertices of opposite sign and two singularities of opposite sign from the
foliation.   
\item The graph of singular leaves for $D_i$ contains a vertex of valence
2, type $(b,b)$ and sign $(+,-)$ (see Figure \ref{typebbvertex}). 
Do a $bb$-exchange move as defined in  \cite{bf}.  
This move deletes two vertices of opposite sign and two singularities of 
opposite sign from the foliation.
\end{enumerate} 
\end{theorem}
The theorem leads to a natural question: Is this set of moves ``minimal'',
or can we eliminate one or more of them?  The data collected in this paper helps to
begin to answer that question.

\begin{question} Do there exist examples of embeddable foliated
discs which have no vertices of valence 1, i.e.  which have no 
end-tiles? If not, then Theorem \ref{theorem:simplifying the disc} 
could obviously be simplified by eliminating moves 2 and 3. 
\end{question}


At the time when Theorem \ref{theorem:simplifying
the disc} was proved, we knew, a classical theorem due
to Magnus and Pelluso \cite{mp} that there are no examples with $n\leq 3$.
We also knew, form Morton's work in \cite{m:4b}, that such
examples \emph{exist} when the braid index is 4. 
Discs whose boundaries have braid index 4 will have $N\geq0$ negative vertices
and $4+N$ positive vertices, so we searched.
There are no such examples of
4-braids with $(P,N) = (5,1), (6,2)$ or $(7,3) $. We found 16 examples with $(P,N) = (8,4)$, 
and (to our great surprise) none with $(P,N) = (9,5)$ or
$(10,6)$. The data suggest that there is structure, not yet understood.

\begin{question}
Do there exist examples of embeddable foliated
discs which have no vertices of valence 1 and no vertices of valence 2, type
$(b,b)$, sign $(+,-)$?  
\end{question}

We do not know the answer to this question. All of our 4-braid examples which lack
vertices of valence 1 have vertices of valence 2 and both type $(a,b)$  and type $(b,b)$
with sign $(+,-)$. 

\begin{figure}
\begin{center}
\includegraphics{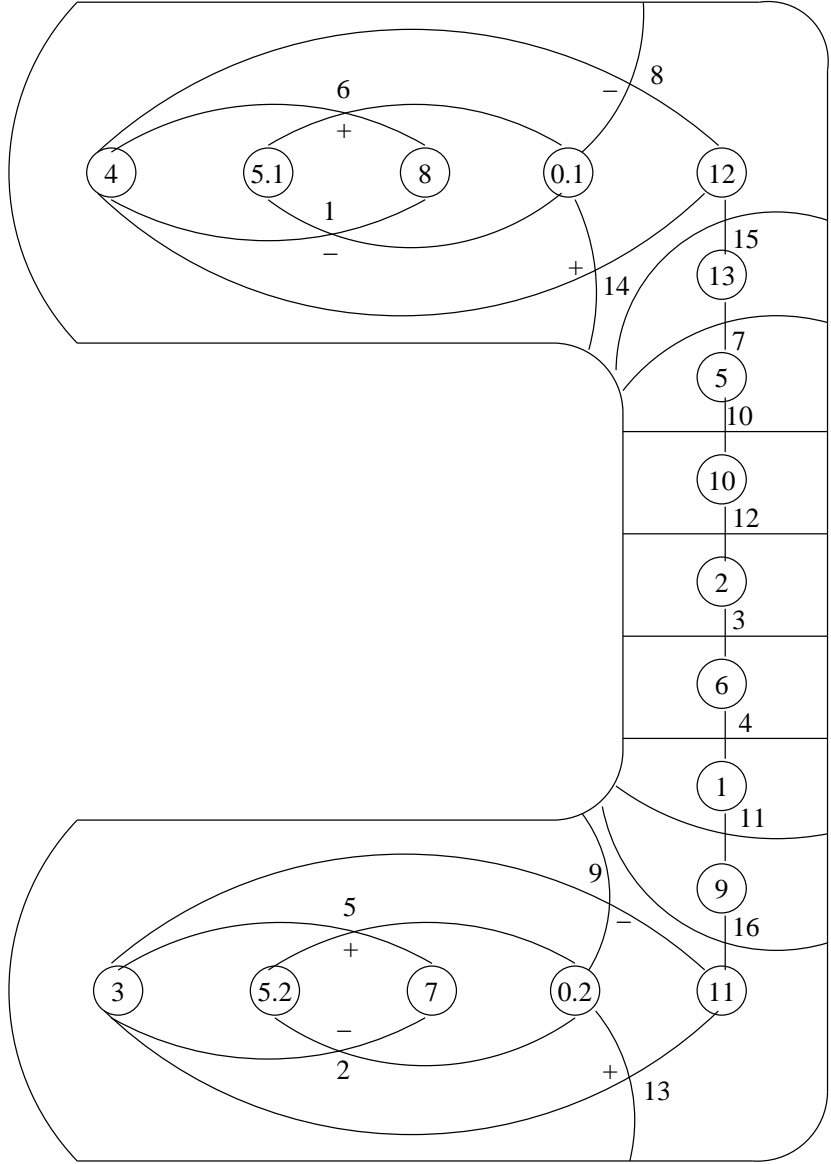}
\end{center}
\caption{An embeddable disc with no vertices of type $(a)$ and no vertices 
of type $(a,b)$ and sign $(+,-)$.}
\label{discodijoan}
\end{figure}

\begin{question}
Do there exist examples of embeddable foliated
discs which have no vertices of valence 1 and no vertices of valence 2, type
$(a,b)$, sign $(+,-)$?  Our final example says ``yes'' (see Figure
\ref{discodijoan}).  To find it we had  to go to $(P,N)=(13,4)$.  This example tells us
that we cannot eliminate  move 3. 
\end{question}

We do not know whether simpler examples of the same type exist.  There are hints of much
more structure in the first set of examples, but we don't have enough data to say more
at this time.  


\begin{question}
Can move 3 be eliminated in the special case $n=4$, i.e. $(P,N)=(N+4,N)$? (Notice that the 
counterexample in Figure \ref{discodijoan} has braid index 9).
\end{question}

\textbf{The data}
There is only one vertex string for $(P,N)=(5,1)$.
Our GAP procedure ComputeCycles(5,1) gives for this 64 cycles, 
none of which is without end-tiles. 

There are four vertex strings for $(P,N)=(6,2)$.
Our GAP procedure ComputeCycles(6,2) gives for them 276 cycles, 
none of which is without end-tiles. 

There are 12 vertex strings with $(7,3)$ vertices, but 4 of them have no
associated cycles. For the other 8 strings, our GAP procedure 
ComputeCycles(7,3) gives 828 cycles, none of which is without end-tiles.

\

The first interesting discs for $B_4$ are those with $(8,4)$ vertices: 
there are 43 vertex strings, 29 of which have no associated cycles; for
the other 14 vertex strings, our GAP procedure 
ComputeCycles(8,4) gives 2944 cycles, 16 of which are with no
end-tiles: all these 16 have the same vertex string (which also has other
12 discs with end-tiles).
The 16 discs without end-tiles are the following 8 and their
inverses:
\be
\item $D_1=\{\ol{[3.1, 5, 8]},\ol{[4, 5.1, 7]},[3.1, 6, 8], 
    \ol{[0.1, 3, 7]},\ol{[0.2, 2, 3.1, 6]}, 
     [1, 4, 5.1],$\\ $[2, 3.1, 5],\ol{[1, 2, 5.1]},
    \ol{[2, 3, 5.1]},[0.2, 2, 6],[0.1, 3, 5.1, 7 ]\},$ 
\item $D_2=\{\ol{[3.1, 5, 8]},\ol{[4, 5.1, 7]},[3.1, 6, 8], 
    \ol{[0.1, 3, 7]},\ol{[0.2, 2, 3.1, 6]}, 
     [1, 4, 5.1],$\\ $[2, 3.1, 5],\ol{[1, 3, 5.1]},
      [ 1, 2 ],[0.2, 2, 6],[0.1, 3, 5.1, 7 ]\},$ 
\item $D_3=\{\ol{[3.1, 5, 8]},\ol{[4, 5.1, 7]},[3.1, 6, 8], 
     \ol{[0.1, 3, 7]},\ol{[0.2, 2, 3.1, 6]}, 
     [1, 4, 5.1],$\\ $[2, 3.1, 5], [ 2, 3 ], 
     \ol{[1, 3, 5.1]},[0.2, 2, 6], 
     [0.1, 3, 5.1, 7 ]\},$ 
\item $D_4=\{\ol{[3.1, 5, 8]},\ol{[4, 5.1, 7]},[3.1, 6, 8], 
     \ol{[0.1, 3, 7]},\ol{[0.2, 2, 3.1, 6]}, 
     [1, 4, 5.1],$\\  $[ 4, 5 ],[2, 3.1, 5], 
     \ol{[1, 3, 5.1]},[0.2, 2, 6], 
     [0.1, 3, 5.1, 7 ]\},$ 
\item $D_5=\{\ol{[3.1, 5, 8]},\ol{[4, 5.1, 7]},[3.1, 7, 8], 
     [3.1, 6, 7],\ol{[0.1, 3, 7]}, 
     \ol{[0.2, 2, 3.1, 6]},$\\ $[1, 4, 5.1], 
     [2, 3.1, 5],\ol{[1, 3, 5.1]},[0.2, 2, 6],
     [0.1, 3, 5.1, 7 ]\},$ 
\item $D_6=\{\ol{[3.1, 5, 8]},\ol{[4, 5.1, 7]},[3.1, 6, 8], 
      [ 7, 8 ],\ol{[0.1, 3, 7]},\ol{[0.2, 2, 3.1, 6]},$\\ 
     $[1, 4, 5.1],[2, 3.1, 5],\ol{[1, 3, 5.1]},
     [0.2, 2, 6],[0.1, 3, 5.1, 7 ]\},$ 
\item $D_7=\{\ol{[3.1, 5, 8]},\ol{[4, 5.1, 7]}, [ 6, 7 ], 
     [3.1, 6, 8],\ol{[0.1, 3, 7]}, 
     \ol{[0.2, 2, 3.1, 6]},$\\ $[1, 4, 5.1], 
     [2, 3.1, 5],\ol{[1, 3, 5.1]},[0.2, 2, 6],
     [0.1, 3, 5.1, 7 ]\},$ 
\item $D_{8}=\{\ol{[3.1, 5, 8]}, [ 4, 5 ],\ol{[4, 5.1, 7]}, 
     [3.1, 6, 8],\ol{[0.1, 3, 7]}, 
     \ol{[0.2, 2, 3.1, 6]},$\\ $[1, 4, 5.1], 
     [2, 3.1, 5],\ol{[1, 3, 5.1]},[0.2, 2, 6],
     [0.1, 3, 5.1, 7 ]\}.$
\ee 
Their corresponding boundary words, depending on the sign assigned
to the $aa$-saddle, are the following:
\be
\item $BW_1=(4,3)\ol{(3,2)}\ol{(4,3)}\ol{(4,3)}
  (3,2)\ol{(2,1)}(3,2)(2,1)\ol{(3,2)},$
\item $BW_2^-=BW_1, (EW_2^-=r_8(EW_1)),$
\item $BW_2^+=(4,3)\ol{(3,2)}\ol{(4,3)}\ol{(4,3)}
  (3,2)\ol{(2,1)}(3,2)(2,1)
  (2,1)(2,1)\ol{(3,2)},$
\item $BW_3^-=BW_1, (EW_3^-=r_8(EW_1)),$
\item $BW_3^+=BW_2^+, (EW_3^+=r_8(EW_2^+)),$
\item $BW_4^-=(4,3)\ol{(3,2)}\ol{(4,3)}\ol{(4,3)}
  (3,2)(2,1)\ol{(3,2)},$
\item $BW_4^+=(4,3)\ol{(3,2)}\ol{(4,3)}\ol{(4,3)}
  (3,2)\ol{(2,1)}(3,2)(3,2)
  (2,1)(2,1)\ol{(3,2)},$
\item $BW_5=(4,3)\ol{(3,2)}\ol{(4,3)}(3,2)
  \ol{(2,1)}(3,2)(2,1)(2,1)
  \ol{(3,2)},$
\item $BW_6^-=(4,3)\ol{(3,2)}\ol{(4,3)}\ol{(4,3)}
  \ol{(4,3)}(3,2)\ol{(2,1)}(3,2)
  (2,1)(2,1)\ol{(3,2)},$
\item $BW_6^+=BW_5, (EW_6^+=r_3(EW_5)),$
\item $BW_7^-=BW_6^-, (EW_7^-=r_3(EW_6^-)),$
\item $BW_7^+=BW_5, (EW_7^+=r_3(EW_5)),$
\item $BW_{8}^-=(4,3)\ol{(3,2)}\ol{(3,2)}\ol{(4,3)}
  \ol{(4,3)}(3,2)\ol{(2,1)}(3,2)
  (2,1)(2,1)\ol{(3,2)},$
\item $BW_{8}^+=\ol{(4,3)}(3,2)\ol{(2,1)}(3,2)
  (2,1)(2,1)\ol{(3,2)}.$
\ee
The conjugacy classes have been computed by S. J. Lee:
considering the 16 words (the 8 which are different in this list, and their inverses),
almost all of them are non conjugate. There are only two pairs of conjugate braids:
$BW_4^-\sim(BW_8^+)^{-1}$ and their inverses.

\begin{figure}
\begin{center}
\includegraphics{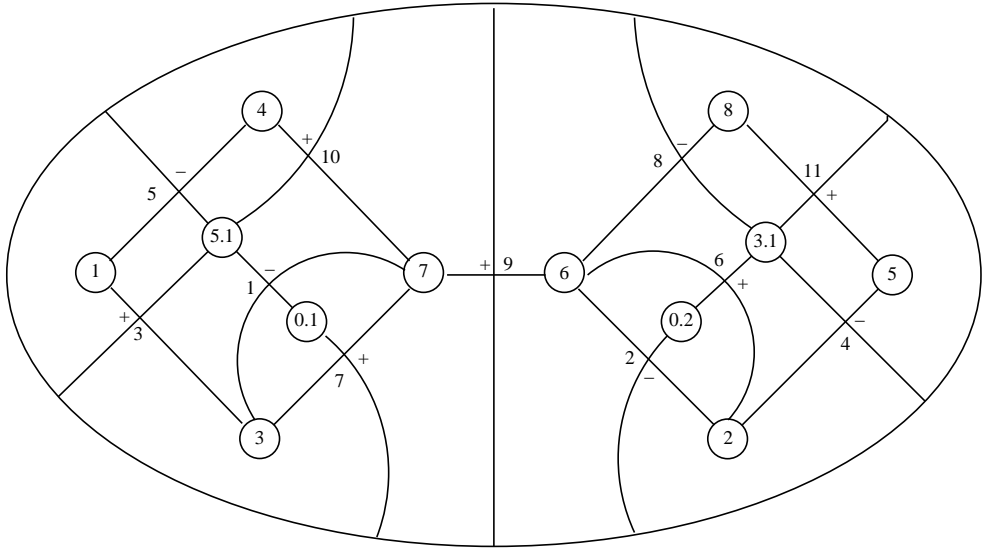}
\end{center}
\caption{The disc $(D_7^-)^{-1}$.}
\label{discogretchen}
\end{figure}

The disc $(D_7^-)^{-1}$, that can be seen in Figure \ref{discogretchen}, is
the disc corresponding to Morton's braid, as shown by G. Wright in \cite{w}.

\begin{figure}
\begin{center}
\includegraphics{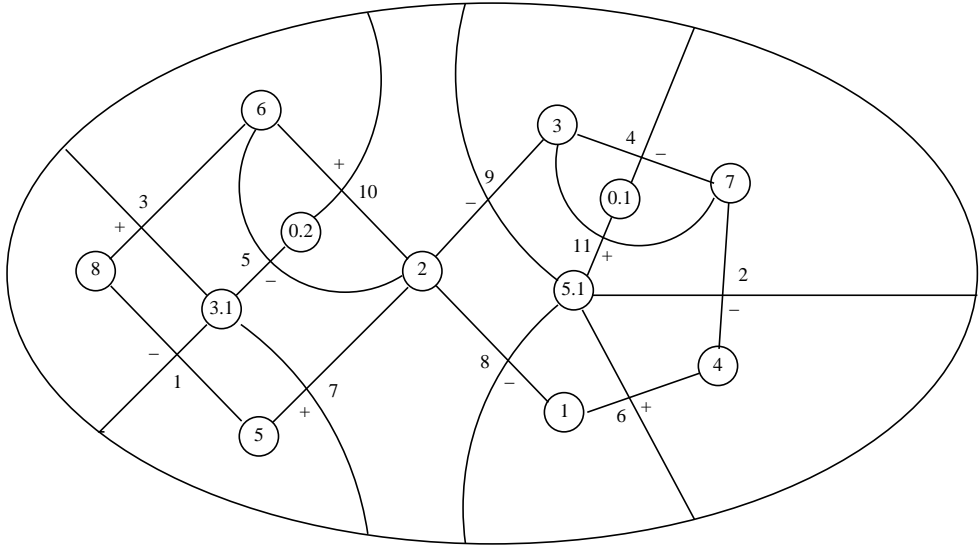}
\end{center}
\caption{The disc $D_1$.}
\label{disco1}
\end{figure}

All other discs have a similar structure, with two `squares' of vertices 
joined by an $aa$-saddle; only discs $D_1,D_5$ and their inverses have
a different structure: a pentagon joined directly to a square, as can be
seen in Figure \ref{disco1}.

\

There are 9,288 cycles with $(9,5)$ vertices, but none of them is 
without end-tiles.
There are 37,952 cycles with $(10,6)$ vertices, but none of them is 
without end-tiles.

\

\begin{figure}
\begin{center}
\includegraphics{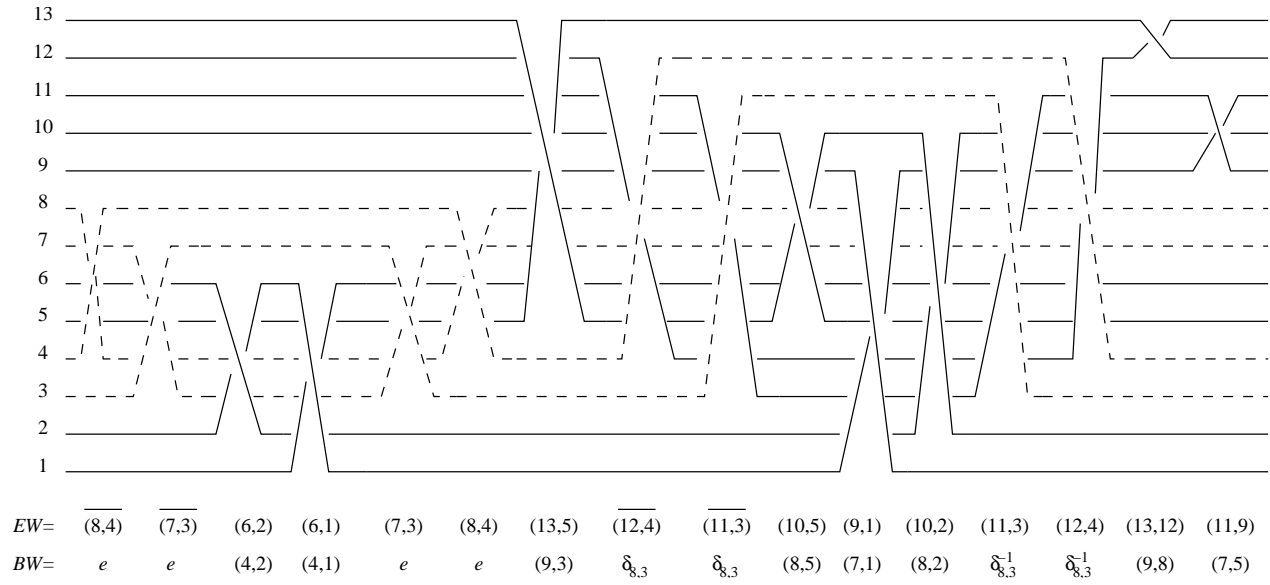}
\end{center}
\caption{The extended braid $EW$; the dashed strands have to be
removed to get the boundary braid $BW$.}
\label{braidjoan}
\end{figure}

An interesting disc with $(13,4)$ vertices, found by Birman and Menasco, is the following
(see Figure \ref{discodijoan}):
$$D=[[[0.1,4,5.1,8],-1],[[3,0.2,7,5.2],-1],[2,6],[1,6],[[3,0.2,7,5.2],1],$$
$$[[0.1,4,5.1,8],1],[5,13],[[0.1,4,12],-1],[[3,0.2,11],-1],[5,10],[1,9],$$
$$[2,10],[[3,0.2,11],1],[[0.1,4,12],1],[12,13],[9,11]].$$
$P=13, N=4, n=9, P+N-1=16$, number of saddles.
$$EW=\ol{(8,4)(7,3)}(6,2)(6,1)(7,3)(8,4)(13,5)\ol{(12,4)(11,3)}(10,5)(9,1)$$
$$(10,2)(11,3)(12,4)(13,12)(11,9).$$
$$\rho(W)=(1,6,10,5,12,13,2,11,9)(3)(4)(7)(8).$$
Deleting the four strands corresponding to the four 1-cycles of
$\rho(W)$, as can be seen in Figure \ref{braidjoan}, get the 9-braid
$$BW=(4,2)(4,1)(9,3)\delta_{8,3}\delta_{8,3}(8,5)(7,1)
(8,2)\delta^{-1}_{8,3}\delta^{-1}_{8,3}(9,8)(7,5)=$$
$$=(4,2)(4,1)(9,3)(8,7)(7,6)(6,5)(5,4)(4,3)(8,7)(7,6)(6,5)(5,4)(4,3)(8,5)$$
$$(7,1)(8,2)\ol{(4,3)(5,4)(6,5)(7,6)(8,7)(4,3)(5,4)(6,5)(7,6)(8,7)}(9,8)(7,5).$$

\section*{Appendix: Commented List of Main GAP Procedures and Functions}

\begin{GAPproc}{EnumeratePositiveGoodWords(p)}%
Enumerates the positive good words (braids) for \T p vertices up to
conjugations by $\delta$ and by subwords. The number of such words is
inherently exponential in \T p, so even if the test for easy conjugations, that
is, the test for equivalence modulo permutation of indices and of letters is
easy (quadratic in \T p), the overall complexity of the procedure is not
polynomial in \T p.
\end{GAPproc}

\begin{GAPproc}{EnumerateGoodWords(p)}%
Enumerates the good words (braids) by suitably inserting signs in the words
obtained through \T{EnumeratePositiveGoodWords(p)}. Each word is a list of
pairs; each pair is formed by the generator name (a pair of increasing indices)
and a sign.
\end{GAPproc}

\begin{GAPproc}{GenerateDiscBoundary(n, w)}%
Takes a good word \T w with \T n indices as input and gives back the
boundary-point list of the corresponding tiled disc. It also checks whether
the word is really a good one, and returns fail if this is not true.
\end{GAPproc}

\begin{GAPproc}{GetInsertionArcs(n, b)}%
Takes the braid index and the list of boundary points of an extended word and
gives back a list of insertion intervals as records with fields \T i, \T f, and
\T I, where \T i is the initial point of the insertion arc, \T f is the final point (they may
coincide) and \T I is the related list of intervals of possible positions for a
new negative vertex.

A disc boundary point is described as a record with fields \T s and \T l, where
\T s describes the saddle associated with the boundary point and \T l is its
level; the saddle is described by a pair whose first coordinate is the list of
vertices involved in the saddle (starting from a positive vertex, and with, if
possible, alternating signs), and the second one is the sign.
\end{GAPproc}

\begin{GAPproc}{GetSaddles(n, s, b)}%
Takes the braid index \T n of an extended word, the number \T s of saddles of
the disc generating the extended word and the list \T b of boundary points (in
the same format of \T{GetInsertionArcs()}, and gives back a list of records
with the following fields: \T i and \T f are the initial and final positive
vertices of the insertion arc; \T N is a list of cyclic intervals, indicating
the possible position of the new negative vertex; \T S is a list of cyclic
intervals, indicating the possible levels of the new saddle; the sign \T s is
the possible sign of the new saddle, \T I is a subinterval of \T N where the
saddle can exist with sign \T s, \T{pi} and \T{pf} are initial and final points
(given by their order on the boundary) of the insertion arc.
\end{GAPproc}

\begin{GAPproc}{Embeddable(n, nsaddles, boundary, bbsaddles)}%
Takes the number \T n of positive vertices of a disc, the number \T{nsaddles} of saddles, 
the list \T{boundary} of boundary points and the list \T{bbsaddles} of $bb$-saddles, 
and outputs whether or not this tiled disc is embeddable. \T{bbsaddles} has the
same format as \T{boundary}, but the order is not relevant.
\end{GAPproc}

\begin{GAPproc}{InsertVertices(n, nsaddles, boundary, bbsaddles, depth)}%
Takes the number \T n of positive vertices of a disc, the number \T{nsaddles} of saddles, 
the list \T{boundary} of boundary points, list \T{bbsaddles} of $bb$-saddles, a
limit \T{depth} on the number of negative vertices to add (if equal to -1, it tries all
possible insertions) and returns a list of embeddable discs with additional
negative vertices; each disc is a record with fields \T b, containing the list
of boundary points, and \T{bb}, containing the list of $bb$-saddles.
\end{GAPproc}

\begin{GAPproc}{EnumerateNecklaces(n, k, d)}%
Returns the set of necklaces of length \T n, \T k colours and density (number
of nonzeros) \T d using the Ruskey--Sawada algorithm~\cite{rs}.  

Even if the number of necklaces is superpolynomial in their length, the
Ruskey--Sawada algorithm has constant amortized complexity, that is, the time
required to generate a single necklace is constant, which implies that the time
necessary to enumerate all necklaces is linear in their number.
\end{GAPproc}

\begin{figure}
\begin{center}
\includegraphics{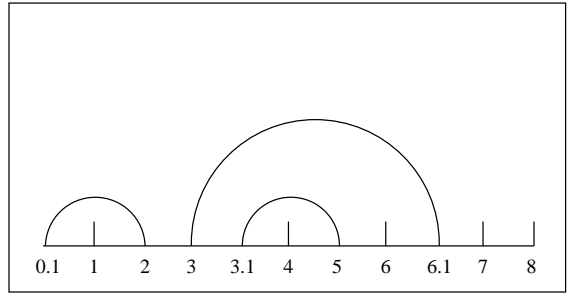}
\end{center}
\caption{An embeddable essential halfplane.}
\label{halfplane}
\end{figure}

\begin{GAPproc}{EnumerateHalfPlanes(vertexString)}%
Generates all possible embeddable essential halfplanes on the given string of
vertices, which must be a vector of zeroes (negative vertices) and ones
(positive vertices). Returns a record with two fields, \T{halfPlane}, described
below, and \T{vertexString}, which copies the input field.

The field \T{halfPlane} is a list of vectors indexed by the vertices. The
$i$-th component of a vector is a record with fields \T v and \T c, where \T v
is the vertex linked by an arc to vertex $i$ and \T c the index of the
``outer'' connected component \T v belongs to (the vertices are numbered all
from $1$ to $P+N$; $0$ means no connection with any vertex).  For instance \T{[
rec( v := 3, c := 0 ), rec( v := 0, c := 1 ), rec(v := 1, c := 0 ), rec( v :=
9, c := 0 ), rec( v := 7, c := 2 ), rec( v := 0, c := 3 ), rec( v := 5, c := 2
), rec( v := 0, c := 2 ), rec(v := 4, c := 0 ), rec( v := 0, c := 0 ), rec( v
:= 0, c := 0 ) ]} describes the halfplane of Figure~\ref{halfplane}.
\end{GAPproc}

\begin{GAPproc}{MakeGraph(halfPlaneRecord)}%
Builds the halfplane adjacency graph for the specified set of halfplanes. 
The result is a record with fields \T{halfPlane}, \T{vertexString}, \T{graph},
\T N, \T A and \T{saddle}.
\T{graph} is an adjacency list: to each halfplane $i$ we associate a list
of records with fields \T v, \T s and \T a where \T v is the index of the next halfplane,
\T s is the associated saddle index (in the list \T{saddle}) multiplied by its sign,
and \T a is the index, in the \T v-th list, of the symmetric edge.
  
The list \T{saddle} is a list of vectors of variable length (2, 3 or 4) indicating
which vertices are involved in a saddle. The vertices are given by their indices
in \T{halfPlaneRecord.vertexString}.

The fields \T N and \T A hold the number of nodes and arcs of the generated graph.

The graph is generated by applying a simple test to each pair of nodes
(halfplanes); thus, the generation requires time quadratic in the number of
halfplanes.
\end{GAPproc}

\begin{GAPproc}{MakeRandomGraph(halfPlaneRecord, isRandom)}%
Same as \T{MakeGraph()}, but \T{isRandom} is a boolean indicating whether we
want to scramble randomly the halfplane order. This is useful for single-cycle
generation.

If \T{isRandom} is false the halfplanes are used in the order in which they
appear in \T{halfPlaneRecord}. Otherwise, they are permuted randomly. The GAP
pseudorandomness seed values \T{R\_N} and \T{R\_X} are accumulated into the two
additional fields of the resulting record with the same name.
\end{GAPproc}

\begin{GAPproc}{EnumerateCycles(halfPlaneGraph, stopAt)}%
This function is the core of the enumeration process; it enumerates the 
$H_\theta$-sequences related to the provided \T{halfPlaneGraph}; the sequences
are expressed as lists of records with fields \T v and \T s, where \T v is a
node and \T s is the signed index of the saddle
that labels the edge towards the next node of the cycle (however, $aa$-saddles
have always positive sign).

The parameter \T{stopAt} specifies how many cycles to generate. If it is zero,
all cycles will be generated. Otherwise, only \T{stopAt} cycles will be
generated. By applying this function with \T{stopAt}=1 on randomly scrambled
graphs, it is possible to generate longer cycles than those allowed by
exhaustive enumeration.

The algorithm used in this function is based on the notion of \textbf{dual
graph}: given a directed graph $G$ with node set $V$ and directed-edge set $A$,
the dual graph has as node set $A$ and directed-edge set $B\subseteq A\times
A$, where $(a,a')\in B$ iff the target of $a$ is the source of $a'$.  More
precisely, if we call $s:A\to V$ and $t:A\to V$ the (obvious) \textbf{source}
and \textbf{target} functions of $G$, we can build the pullback
\[
\begin{array}{ccc}
  B & \longrightarrow & A \\
  \downarrow & & \downarrow\mbox{\scriptsize\textit{s}}\\
  A & \stackrel t\longrightarrow & V
\end{array}
\]

Then, the dual graph has directed-edge set $B$, and the projections of the
pullback are precisely its source and target functions. The crucial observation
is that \emph{elementary} cycles in the dual graph (i.e., cycles that never
pass twice through the same \emph{node}) are in bijection with cycles of the original
graph that never pass twice through the same \emph{edge}. This allows one to
use standard enumeration methods for elementary cycles for enumerating
$H_\theta$-sequences. 

Nevertheless, sophisticated methods such as Johnson's algorithm~\cite{jo} turn
out to be inefficient in our case: this happens because the length of the
cycles to be generated is fixed, and is very small with respect to the size of
the graph. After several experiments, we focused on a two-phase 
visiting algorithm. In the first phase, we use standard depth-first
enumeration techniques to generate candidate loop-free cycles of length shorter
than or equal to $P+N-1$; in the second phase, we enrich each candidate with
loops, and test the various conditions that must be satisfied to obtain an
$H_\theta$-sequence.

In the first phase, a total ordering is established on the nodes of the dual
graph. At each round of the generation, a \emph{root} node of the dual graph is
chosen as first node of the cycle to be generated; moreover, the generation
process only considers nodes larger than the root.  This ensures that loop-free
cycles are never generated twice. Moreover, a precomputation of the distances
from each node to the root (using a standard breadth-first visit) allows to cut
prematurely cycles that could never ``get back in time'' because they have
moved too far apart from the root.

In the second phase (which is invoked for each candidate loop-free cycle) we
exhaustively try to add loops so to obtain an $H_\theta$-sequence.

Note that the number of cycles to be generated is inherently superexponential;
it is also very difficult to estimate the amortized complexity of the
enumeration process. The division in two phases cuts a large part of the search
space with respect to a simple enumeration, but nonetheless exhaustive
enumeration is possible only for a very small number of vertices.
\end{GAPproc}

  

\begin{GAPproc}{ComputeCycles(p, n)}%
This function generates all cycles for \T p positive and \T n negative
vertices, using the functions above. The result is provided in two fields
named \T c and \T e, which contains all cycles, and all cycles that do not
contain end tiles, respectively. A cycle here is specified in a more
user-friendly form, that is, as a sequence of saddles, each saddle being
of the form \T{[v,w]} for $aa$-saddles, and of the form
\T{[[v,w,z],s]} or \T{[[u,v,w,z],s]} for $ab$ and $bb$-saddles (\T s is the
sign of the saddle). Vertices are named as in this paper.
\end{GAPproc}

\begin{GAPproc}{PSFilm(c, halfPlaneGraph)}%
Outputs in the current directory a file named \T{film.ps} containing a
PostScript visualization of the sequence of halfplanes traversed by the
$H_\theta$-sequences represented by the cycle \T c.  The visualization also
contains one of the corresponding boundary braid words ($aa$-saddles signs are
arbitrary).
\end{GAPproc}

\begin{GAPproc}{DrawDisc(c, halfPlaneGraph)}%
Outputs in the current directory a file named \T{disc.gml} containing the graph
structure of the tiled disc associated to the cycle \T c. Note that the graph
is \emph{not yet} embedded in the plane: a planar embedding layout algorithm
must be applied.
\end{GAPproc}

\noindent Joan S. Birman, Dept.of Mathematics, Barnard College\\
Mail Code 4427, Columbia University\\ 2990 Broadway, New York, N.Y., USA 10027 \\
e-mail: jb@math.columbia.edu\\

\noindent Marta Rampichini,
Dipartimento di Matematica\\Universit\`a di
Milano\\via Saldini 50, 20133, Milano, Italy\\
e-mail: rampichini@mat.unimi.it\\

\noindent Paolo Boldi and Sebastiano Vigna, Dipartimento di Scienze
dell'Informazione\\
Universit\`a di Milano\\via Comelico 39/41, 20135 Milano, Italy\\
e-mail: \{boldi,vigna\}@dsi.unimi.it
\end{document}